# TESTING FOR JUMPS IN A DISCRETELY OBSERVED PROCESS

By Yacine Aït-Sahalia[1] and Jean Jacod

*Princeton University and Université Pierre et Marie Curie*

We propose a new test to determine whether jumps are present in asset returns or other discretely sampled processes. As the sampling interval tends to 0, our test statistic converges to 1 if there are jumps, and to another deterministic and known value (such as 2) if there are no jumps. The test is valid for all Itô semimartingales, depends neither on the law of the process nor on the coefficients of the equation which it solves, does not require a preliminary estimation of these coefficients, and when there are jumps the test is applicable whether jumps have finite or infinite-activity and for an arbitrary Blumenthal–Getoor index. We finally implement the test on simulations and asset returns data.

**1. Introduction.** The problem of deciding whether the continuous-time process which models an economic or financial time series should have continuous paths or exhibit jumps is becoming an increasingly important issue, in view of the high-frequency observations that are now widely available. In the case where a large jump occurs, a simple glance at the dataset might be sufficient to decide this issue. But such large jumps are usually infrequent, may not belong to the model itself, can be considered as breakdowns in the homogeneity of the model, or may be dealt with separately using other methods such as risk management.

On the other hand, a visual inspection of most such time series in practice does not provide clear evidence for either the presence or the absence of small or medium sized jumps. Since small frequent jumps should definitely be incorporated into the model, and since models with and without jumps do have quite different mathematical properties and financial consequences (for option hedging, portfolio optimization, etc.), it is important to have statistical methods that can shed some light on the issue.

Received October 2006; revised October 2007.
[1]Supported in part by NSF Grants SBR-03-50772 and DMS-05-32370.

Determining whether a process has jumps has been considered by a number of authors. Aït-Sahalia (2002) relies on restrictions on the transition function of the process that are compatible with continuity of the process, or lack thereof, to derive a test for the presence of jumps at any observable frequency. Using high-frequency data, Carr and Wu (2003) exploit the differential behavior of short dated options to test for the presence of jumps. Multipower variations can separate the continuous part of the quadratic variation; see Barndorff-Nielsen and Shephard (2006). Andersen, Bollerslev and Diebold (2003) and Huang and Tauchen (2006) study financial datasets using multipower variations, in order to assess the proportion of quadratic variation attributable to jumps. Jiang and Oomen (2005) construct a test motivated by the hedging error of a variance swap replication strategy. Other methods have been introduced as well [see, e.g., Lee and Mykland (2008)] and the literature about evaluating the volatility when there are jumps [see, e.g., Aït-Sahalia (2004), Aït-Sahalia and Jacod (2007), Mancini (2001), Mancini (2004) and Woerner (2006b)] can also be viewed as an indirect way of checking for jumps. Woerner (2006a) proposes estimators of the Blumenthal–Getoor index or the Hurst exponent of a stochastic process based on high-frequency data, and those can be used to detect jumps and their intensity.

These are two closely related but different issues: one is to decide whether jumps are present or not, another one is to determine the impact of jumps on the overall variability of the observed process. So far, most of the literature has concentrated on the second issue, and it usually assumes a special structure on the process, especially regarding its jump part if one is present, such as being a compound Poisson, or sometimes a Lévy process. We plan to have a systematic look at the question of measuring the impact of jumps in a future paper, using a variety of distances and doing so under weak assumptions, and also taking into account the usually nonnegligible market microstructure noise.

In this paper, we concentrate on the first problem, namely whether there are jumps or not, while pretending that the underlying process is perfectly observed at $n$ discrete times, and $n$ is large. We introduce a direct and very simple test which gives a solution to this problem, irrespective of the precise structure of the process within the very large class of Itô semimartingales.

More specifically, a process $X = (X_t)$ on a given time interval $[0, T]$ is observed at times $i\Delta_n$ for $\Delta_n = T/n$. We propose an easy-to-compute family of test statistics, say $S_n$, which converge as $\Delta_n \to 0$ to 1 if there are jumps, and to another deterministic and known value (such as 2) if there are no jumps. This holds as soon as the process $X$ is an Itô semimartingale, and it depends neither on the law of the process nor on the coefficients of the equation which it solves, and it does not require any preliminary estimation of these coefficients. We provide a central limit theorem for $S_n$, under both alternatives (jumps and no jumps); again this does not require any a priori knowledge of the coefficients of the model, and when there are jumps it is applicable whether the jumps have finite or infinite activity, and for an arbitrary Blumenthal–Getoor index. Hence we can construct tests with a given level of significance asymptotically and which are fully nonparametric.

The paper is organized as follows. Sections 2 and 3 describe our setup and the statistical problem. We provide central limit theorems for our proposed statistics in Section 4, and use them to construct the actual tests in Section 5. We present the results of Monte Carlo simulations in Section 6, and in Section 7 examine the empirical distribution of the test statistic over all 2005 transactions of the Dow Jones stocks. Proofs are in Section 8.

**2. Setting and assumptions.** The underlying process $X$ which we observe at discrete times is a one-dimensional process which we specify below. Observe that taking a one-dimensional process is not a restriction in our context since, if it were multidimensional, a jump would necessarily be a jump for at least one of its components, so any test for jumps can be performed separately on each of its components.

As already mentioned, we do not want to make any specific model assumption on $X$, such as assuming some parametric family of models. We do need, however, a mild structural assumption which is satisfied in all continuous-time models used in finance, at least as long as one wants to rule out arbitrage opportunities. In any case, in the absence of some kind of assumption, anything can happen from a continuous process which is a linear interpolation between the observations, to a pure jump process which is constant between successive observation times.

Our structural assumption is that $X$ is an Itô semimartingale on some filtered space $(\Omega, \mathcal{F}, (\mathcal{F}_t)_{t\geq 0}, \mathbb{P})$, which means it can be written as

$$(1) \quad \begin{aligned} X_t = X_0 &+ \int_0^t b_s\, ds + \int_0^t \sigma_s\, dW_s + \int_0^t \int_E \kappa \circ \delta(s,x)(\mu - \nu)(ds, dx) \\ &+ \int_0^t \int_E \kappa' \circ \delta(s,x) \mu(ds, dx), \end{aligned}$$

where $W$ and $\mu$ are a Wiener process and a Poisson random measure on $\mathbb{R}_+ \times E$ with $(E, \mathcal{E})$ an auxiliary measurable space on the space $(\Omega, \mathcal{F}, (\mathcal{F}_t)_{t\geq 0}, \mathbb{P})$ and the predictable compensator (or intensity measure) of $\mu$ is $\nu(ds, dx) = ds \otimes \lambda(dx)$ for some given finite or $\sigma$-finite measure on $(E, \mathcal{E})$. Moreover $\kappa$ is a continuous function with compact support and $\kappa(x) = x$ on a neighborhood of 0, and $\kappa'(x) = x - \kappa(x)$.

Above, the function $\kappa$ is arbitrary; a change of this function amounts to a change in the drift coefficient $b_t$. Also, it is always possible to take $(E, \mathcal{E}, \lambda)$ to be $\mathbb{R}$ equipped with Lebesgue measure.

Of course the coefficients $b_t(\omega)$, $\sigma_t(\omega)$ and $\delta(\omega, t, x)$ should be such that the various integrals in (1) make sense [see, e.g., Jacod and Shiryaev (2003) for a precise definition of the last two integrals in (1)], and in particular $b_t$ and $\sigma_t$ are optional processes and $\delta$ is a predictable function. However, we need a bit more than just the minimal integrability assumptions, plus the nontrivial fact that the volatility $\sigma_t$ is also an Itô semimartingale, of the form

$$
\begin{aligned}
\sigma_t = \sigma_0 &+ \int_0^t \widetilde{b}_s \, ds + \int_0^t \widetilde{\sigma}_s \, dW_s + \int_0^t \widetilde{\sigma}'_s \, dW'_s \\
&+ \int_0^t \int_E \kappa \circ \widetilde{\delta}(s,x)(\mu - \nu)(ds, dx) \\
&+ \int_0^t \int_E \kappa' \circ \widetilde{\delta}(s,x)\mu(ds, dx),
\end{aligned}
\tag{2}
$$

where $W'$ is another Wiener process independent of $(W, \mu)$.

All our additional requirements are expressed in the next assumption, for which we need a few additional notations. We write

$$\Delta X_s = X_s - X_{s-} \tag{3}$$

for the jumps of the $X$ process. Of course, $\Delta X_s = 0$ for all $s$ when $X$ is continuous, and for all $s$ outside a (random) countable set in all cases. Let

$$
\delta'_t(\omega) = \begin{cases} \int \kappa \circ \delta(\omega, t, x)\lambda(dx), & \text{if the integral makes sense,} \\ +\infty, & \text{otherwise} \end{cases}
\tag{4}
$$

and $\tau = \inf(t; \Delta X_t \neq 0)$. In the finite-activity case for jumps we have $\tau > 0$ a.s., whereas in the infinite-activity case we usually (but not necessarily) have $\tau = 0$ a.s. Observe that the fifth term in (2) is a.s. equal to $-\int_0^t \delta'_s \, ds$ for all $t$ in the (possibly empty) set $[0, \tau]$; so outside a null set in $(\omega, t)$ for the product measure $\mathbb{P}(d\omega) \otimes dt$, the process $\delta'_t(\omega)$ takes finite values for all $t \leq \tau(\omega)$, whereas we may very well have $\delta'_t(\omega) = \infty$ for all $t > \tau(\omega)$.

ASSUMPTION 1.　(a) All paths $t \mapsto \widetilde{b}_t(\omega)$ are locally bounded.

(b) All paths $t \mapsto b_t(\omega)$, $t \mapsto \widetilde{\sigma}_t(\omega)$, $t \mapsto \widetilde{\sigma}'_t(\omega)$ are right-continuous with left limits.

(c) All paths $t \mapsto \delta(\omega, t, x)$ and $t \mapsto \widetilde{\delta}(\omega, t, x)$ are left-continuous with right limits.

(d) All paths $t \mapsto \sup_{x \in E} \frac{|\delta(\omega, t, x)|}{\gamma(x)}$ and $t \mapsto \sup_{x \in E} \frac{|\widetilde{\delta}(\omega, t, x)|}{\gamma(x)}$ are locally bounded, where $\gamma$ is a (nonrandom) nonnegative function satisfying $\int_E (\gamma(x)^2 \wedge 1)\lambda(dx) < \infty$.

(e) All paths $t \mapsto \delta'_t(\omega)$ are left-continuous with right limits on the semiopen set $[0, \tau(\omega))$.

(f) We have $\int_0^t |\sigma_s| \, ds > 0$ a.s. for all $t > 0$.

The nondegeneracy condition (f) says that almost surely the continuous martingale part of $X$ is not identically 0 on any interval $[0,t]$ with $t>0$. Most results of this paper are true without (f), but not all, and in any case this condition is satisfied in all applications we have in mind. Apart from this nondegeneracy condition, Assumption 1 accommodates virtually all models for stochastic volatility, including those with jumps, and allows for correlation between the volatility and the asset price processes. For example, if we consider a $d$-dimensional equation

$$dY_t = f(Y_{t-})\,dZ_t, \tag{5}$$

where $Z$ is a multidimensional Lévy process with Gaussian components, and $f$ is a $C^2$ function with at most linear growth, then any of the components of $Y$ satisfies Assumption 1 except perhaps (e); this comes from Itô's formula and from the representation of a Lévy process in terms of a Wiener process and a Poisson random measure. The same holds for more general equations driven by a Wiener process and a Poisson random measure.

Except when stated otherwise, we will maintain Assumption 1 throughout the paper and will therefore omit it in the statements of the results that follow.

## 3. The statistical problem.

3.1. *Preliminaries.* Our process $X$ is discretely observed over a given time interval $[0,t]$, and it is convenient in practice to allow $t$ to vary. So we suppose that $X$ can be observed at times $i\Delta_n$ for all $i=0,1,\ldots$ and we take into account only those observation times $i\Delta_n$ smaller than or equal to $t$, whether $t$ is a multiple of $\Delta_n$ or not. Moreover the testing procedures given below are "asymptotic," in the sense that we can specify the level or the power function asymptotically as $n\to\infty$ and $\Delta_n\to 0$.

Recall that we want to decide whether there are jumps or not, for the process (1); equivalently, we want to decide whether the coefficient $\delta$ is identically 0 or not. A few remarks can be stated right away, which show some of the difficulties or peculiarities of this statistical problem:

1. It is a nonparametric problem: we do not specify the coefficients $b$, $\sigma$, $\delta$.
2. It is an asymptotic problem, which only makes sense for high-frequency data.
3. When "$n$ is infinite," that is, in the ideal although unrealistic situation of a complete observation of the path of $X$ over $[0,t]$, we can of course tell whether our particular path has jumps or not. However, when the measure $\lambda$ is finite there is a positive probability that the path $X(\omega)$ has no jump on $[0,t]$, although the model itself may allow for jumps.

4. In the realistic case $n < \infty$ we cannot do better than in the "completely observed case." That is, we can hopefully infer something about the jumps which actually occurred for our particular observed path, but nothing about those which belong to the model but did not occur on the observed path.
5. Rather than "testing for jumps," there are cases in practice where one wants to estimate in some sense the part of the variability of the process which is due to the jumps, compared to the part due to the continuous component. This is what most authors have studied so far, in some special cases at least, and we will take a systematic look at this question in our further paper in preparation on the topic, using essentially the same tools as in the present paper. However, under the null hypothesis that jumps are present, it is not clear how one should go about specifying the proportion of quadratic variation attributable to jumps without already assuming not only the type but also the "quantity" of jumps, for instance when using statistics such as those in Barndorff-Nielsen and Shephard (2006). The test we propose below does not have this problem.
6. Coming back to testing for jumps, an important property of test statistics is that they should be scale-invariant (invariant if $X$ is multiplied by an arbitrary constant). It would also be desirable for the limiting behavior of the statistic to be independent of the dynamics of the process. We will see that our test has all these features.

In view of comments 3 and 4 above, the problem which we really try to solve in this paper is to decide, on the basis of the observations $X_{i\Delta_n}$ which belong to the time interval $[0,t]$, in which of the following two complementary sets the path which we have discretely observed falls:

(6) $\quad \begin{cases} \Omega_t^j = \{\omega : s \mapsto X_s(\omega) \text{ is discontinuous on } [0,t]\}, \\ \Omega_t^c = \{\omega : s \mapsto X_s(\omega) \text{ is continuous on } [0,t]\}. \end{cases}$

If we decide on $\Omega_t^j$, then we also implicitly decide that the model has jumps, whereas if we decide on $\Omega_t^c$ it does not mean that the model is continuous, even on the interval $[0,t]$ (of course, in both cases we can say nothing about what happens after $t$!).

3.2. *Measuring the variability of $X$.* Let us now introduce a number of processes which all measure some kind of variability of $X$, or perhaps its continuous and jump components separately, and depend on the whole (unobserved) path of $X$:

(7) $\qquad A(p)_t = \int_0^t |\sigma_s|^p \, ds, \qquad B(p)_t = \sum_{s \leq t} |\Delta X_s|^p,$

where $p$ is a positive number.

Note that $A(p)$ is finite-valued for all $p > 0$ (under Assumption 1), whereas $B(p)$ is finite-valued if $p \geq 2$ but often not when $p < 2$. Also, recall that $[X, X] = A(2) + B(2)$.

We have $\Omega_t^j = \{B(p)_t > 0\}$ for any $p > 0$, so in a sense our problem "simply" amounts to determining whether $B(p)_t > 0$ for our particular observed path, and with any prespecified $p$. Moreover a reasonable measure of the relative variability, or variance, due to the jumps is $B(2)_t/[X, X]_t$, and this is the measure used, for example, by Huang and Tauchen (2006); other measures of this variability could be $B(p)_t/A(p)_t$ or $B(p)_t/[X, X]_t^{p/2}$ for other values of $p$ (the power $p/2$ in the denominator is to ensure the scale-invariance).

In any event, everything boils down to estimating, on the basis of the actual observations, the quantity $B(p)_t$ in (7), and the difficulty of this estimation depends on the value of $p$. Let

$$(8) \qquad \Delta_i^n X = X_{i\Delta_n} - X_{(i-1)\Delta_n}$$

denote the observed discrete increments of $X$ [all of them, not just those due to jumps, unlike (3)] and define for $p > 0$ the estimator

$$(9) \qquad \widehat{B}(p, \Delta_n)_t := \sum_{i=1}^{[t/\Delta_n]} |\Delta_i^n X|^p.$$

For $r \in (0, \infty)$, let

$$(10) \qquad m_r = \mathbb{E}(|U|^r) = \pi^{-1/2} 2^{r/2} \Gamma\left(\frac{r+1}{2}\right)$$

denote the $r$th absolute moment of a variable $U \sim N(0, 1)$. We have the following convergences in probability, locally uniform in $t$:

$$(11) \qquad \begin{cases} p > 2 & \Rightarrow \quad \widehat{B}(p, \Delta_n)_t \xrightarrow{\mathbb{P}} B(p)_t, \\ p = 2 & \Rightarrow \quad \widehat{B}(p, \Delta_n)_t \xrightarrow{\mathbb{P}} [X, X]_t, \\ p < 2 & \Rightarrow \quad \dfrac{\Delta_n^{1-p/2}}{m_p} \widehat{B}(p, \Delta_n)_t \xrightarrow{\mathbb{P}} A(p)_t, \\ X \text{ is continuous} & \Rightarrow \quad \dfrac{\Delta_n^{1-p/2}}{m_p} \widehat{B}(p, \Delta_n)_t \xrightarrow{\mathbb{P}} A(p)_t. \end{cases}$$

These properties are known: for $p = 2$ this is the convergence of the realized quadratic variation, for $p > 2$ this is due to Lepingle (1976) for all semimartingales, and for the other cases one may see, e.g., Jacod (2008).

The intuition for the behavior of $\widehat{B}(p, \Delta_n)_t$ is as follows. Suppose that $X$ can jump. Among the increments of $X$, there are those that contain a large jump and those that do not. While the increments containing large jumps are much less frequent than those that contain only a Brownian contribution and many small jumps, or only a Brownian contribution when $\lambda$ is finite, they are so much bigger than the rest that when jumps occur, their contribution to $B(p)$ for $p > 2$ overwhelms everything else. This is because high powers ($p > 2$) magnify the large increments at the expense of the small ones. Then the sum behaves like the sum coming from the jumps only; this is the first result in (11). When $p$ is small ($p < 2$), on the other hand, the magnification of the large increments by the power is not strong enough to overcome the fact that there are many more small increments than large ones. Then the behavior of the sum is driven by the summation of all these small increments; this is the third result in (11). When $p = 2$, we are in the situation where these two effects (magnification of the relatively few large increments vs. summation of many small increments) are of the same magnitude; this is the second result in (11). When $X$ is continuous, we are only summing small increments and we get for all values of $p$ the same behavior as the one in the third result where the summation of the small increments dominates the sum; this is the fourth result in (11).

3.3. *The test statistics.* Based on this intuition, upon examining (11), we see that when $p > 2$ the limit of $\widehat{B}(p, \Delta_n)_t$ does not depend on the sequence $\Delta_n$ going to 0, and it is strictly positive if $X$ has jumps between 0 and $t$. On the other hand, when $X$ is continuous on $[0, t]$, then $\widehat{B}(p, \Delta_n)_t$ converges again to a limit not depending on $\Delta_n$, but only after a normalization which does depend on $\Delta_n$.

These considerations lead us to compare $\widehat{B}(p, \Delta_n)_t$ on two different $\Delta_n$-scales. Specifically, we choose an integer $k$ and compare $\widehat{B}(p, \Delta_n)_t$ with $\widehat{B}(p, k\Delta_n)_t$, the latter obtained by considering only the increments of $X$ between successive multiples of $k\Delta_n$. Then we set

$$(12) \qquad \widehat{S}(p, k, \Delta_n)_t = \frac{\widehat{B}(p, k\Delta_n)_t}{\widehat{B}(p, \Delta_n)_t}$$

as our (family of) test statistics.

In view of the first and fourth limits in (11), we readily get:

THEOREM 1. *Let $t > 0$, $p > 2$ and $k \geq 2$. Then the variables $\widehat{S}(p, k, \Delta_n)_t$ converge in probability to the variable $S(p, k)_t$ defined by*

$$(13) \qquad S(p, k)_t = \begin{cases} 1, & \text{on the set } \Omega_t^j, \\ k^{p/2-1}, & \text{on the set } \Omega_t^c. \end{cases}$$

(On the set $\Omega_t^c$, the convergence holds for $p \leq 2$ as well.) Therefore our test statistics will converge to 1 in the presence of jumps and, with the selection of $p = 4$ and $k = 2$, to 2 in the absence of jumps. The following corollary is immediate:

COROLLARY 1. *Let $t > 0$, $p > 2$ and $k \geq 2$. The decision rule defined by*

$$\Im(n) = \Im(n,t,p,k) = \begin{cases} X \text{ is discontinuous on } [0,t], & \text{if } \widehat{S}(p,k,\Delta_n)_t < a, \\ X \text{ is continuous on } [0,t], & \text{if } \widehat{S}(p,k,\Delta_n)_t \geq a \end{cases}$$

*is consistent, in the sense that the probability that it gives the wrong answer tends to $0$ as $\Delta_n \to 0$, for any choice of $a$ in the interval $(1, k^{p/2} - 1)$.*

## 4. Central limit theorems.

4.1. *CLT for power variations.* The previous corollary provides the first step toward constructing a test for the presence or absence of jumps, but it is hardly enough. To construct tests, we need to derive the rates of convergence and the asymptotic variances when $X$ jumps and when it is continuous.

We start with the following general theorem, to be proved in Section 8 (as above, $k$ is an integer larger than 1). The asymptotic variances in the theorem involve the process $A(p)_t$ already defined in (7) as well as the more complex process

$$(14) \qquad D(p)_t = \sum_{s \leq t} |\Delta X_s|^p (\sigma_{s-}^2 + \sigma_s^2)$$

for $p > 0$. Like $B(p)_t$, $D(p)_t$ is finite-valued if $p \geq 2$ but often not when $p < 2$.

THEOREM 2. (a) *Let $p > 3$. For any $t > 0$ the pair of variables*

$$\Delta_n^{-1/2} (\widehat{B}(p,\Delta_n)_t - B(p)_t, \widehat{B}(p,k\Delta_n)_t - B(p)_t)$$

*converges stably in law to a bidimensional variable of the form $(Z(p)_t, Z(p)_t + Z'(p,k)_t)$, where both $Z(p)_t$ and $Z'(p,k)_t$ are defined on an extension $(\widetilde{\Omega}, \widetilde{\mathcal{F}}, (\widetilde{\mathcal{F}}_t)_{t \geq 0}, \widetilde{\mathbb{P}})$ of the original filtered space $(\Omega, \mathcal{F}, (\mathcal{F}_t)_{t \geq 0}, \mathbb{P})$ and conditionally on $\mathcal{F}$ are centered, with $Z'(p)_t$ having the following conditional variance:*

$$(15) \qquad \widetilde{\mathbb{E}}(Z'(p,k)_t^2 \mid \mathcal{F}) = \frac{k-1}{2} p^2 D(2p-2)_t.$$

*Moreover if the processes $\sigma$ and $X$ have no common jumps, the variable $Z'(p)_t$ are $\mathcal{F}$-conditionally Gaussian.*

(b) *Assume in addition that $X$ is continuous, and let $p \geq 2$. The pair of variables*

$$\Delta_n^{-1/2} (\Delta_n^{1-p/2} \widehat{B}(p,\Delta_n)_t - m_p A(p)_t, \Delta_n^{1-p/2} \widehat{B}(p,k\Delta_n)_t - k^{p/2-1} m_p A(p)_t)$$

*converges stably in law to a bidimensional variable $(Y(p)_t, Y'(p,k)_t)$ defined on an extension $(\widetilde{\Omega}, \widetilde{\mathcal{F}}, (\widetilde{\mathcal{F}}_t)_{t \geq 0}, \widetilde{\mathbb{P}})$ of the original filtered space $(\Omega, \mathcal{F}, (\mathcal{F}_t)_{t \geq 0}, \mathbb{P})$* *and which conditionally on $\mathcal{F}$ is a centered Gaussian variable with variance-covariance given by*

$$
(16) \quad \begin{cases} \widetilde{\mathbb{E}}(Y(p)_t^2 \mid \mathcal{F}) = (m_{2p} - m_p^2)A(2p)_t, \\ \widetilde{\mathbb{E}}(Y'(p,k)_t^2 \mid \mathcal{F}) = k^{p-1}(m_{2p} - m_p^2)A(2p)_t, \\ \widetilde{\mathbb{E}}(Y(p)_t Y'(p,k)_t \mid \mathcal{F}) = (m_{k,p} - k^{p/2} m_p^2)A(2p)_t \end{cases}
$$

*and where*

$$
(17) \quad m_{k,p} = \mathbb{E}(|U|^p \, |U + \sqrt{k-1}\,V|^p)
$$

*for $U$, $V$ independent $N(0,1)$ variables.*

The "stable convergence in law" mentioned above is a mode of convergence introduced by Rényi (1963), which is slightly stronger than the mere convergence in law, and its important feature for us is that if $V_n$ is any sequence of variables converging in probability to a limit $V$ on the space $(\Omega, \mathcal{F}, (\mathcal{F}_t)_{t \geq 0}, \mathbb{P})$, whereas the variables $V'_n$ converge stably in law to $V'$, then the pair $(V_n, V'_n)$ converges stably in law again to the pair $(V, V')$.

The property "the processes $\sigma$ and $X$ have no common jumps" may hold, even though both processes are driven by the same Poisson measure; indeed, it holds when the product $\delta \widetilde{\delta}$ vanishes identically, or more generally when $(\delta \widetilde{\delta})(\omega, s, z) = 0$ for $\mathbb{P}(d\omega) \times ds \times \lambda(dz)$ almost all $(\omega, s, z)$. Because of the freedom we have in the choice of the driving Poisson measure, in this case we can also use two independent Poisson measures $\mu$ and $\widetilde{\mu}$ to drive $\sigma$ and $X$, without changing the law of the pair $(X, \sigma)$, whereas conversely if this is true we can "aggregate" the two measures $\mu$ and $\widetilde{\mu}$ into a single one.

The reader will observe the restriction $p > 3$ in (a) above: the CLT simply does not hold if $p \leq 3$, when there are jumps. On the other hand, (b) holds also if $p < 2$, under the additional assumption that $\sigma$ does not vanish, but we do not need this improvement here.

4.2. *CLT for the nonstandardized statistics.* This theorem allows us to deduce a CLT for our statistics $\widehat{S}(p, k, \Delta_n)_t$:

THEOREM 3. (a) *Let $p > 3$ and $t > 0$. Then $\Delta_n^{-1/2}(\widehat{S}(p, k, \Delta_n)_t - 1)$ converges stably in law, in restriction to the set $\Omega_t^j$ of (13), to a variable $S(p, k)_t^j$ which, conditionally on $\mathcal{F}$, is centered with variance*

$$
(18) \quad \widetilde{\mathbb{E}}((S(p,k)_t^j)^2 \mid \mathcal{F}) = \frac{(k-1)p^2}{2} \frac{D(2p-2)_t}{B(p)_t^2}.
$$

*Moreover if the processes $\sigma$ and $X$ have no common jumps, the variable $S(p, k)_t^j$ is $\mathcal{F}$-conditionally Gaussian.*

(b) *Assume in addition that $X$ is continuous, and let $p \geq 2$ and $t > 0$. Then $\Delta_n^{-1/2}(\widehat{S}(p,k,\Delta_n)_t - k^{p/2-1})$ converge stably in law to a variable $S(p,k)_t^c$ which, conditionally on $\mathcal{F}$, is centered normal with variance*

$$(19) \qquad \widetilde{\mathbb{E}}((S(p,k)_t^c)^2 \mid \mathcal{F}) = M(p,k)\frac{A(2p)_t}{A(p)_t^2},$$

*where*

$$(20) \quad M(p,k) = \frac{1}{m_p^2}(k^{p-2}(1+k)m_{2p} + k^{p-2}(k-1)m_p^2 - 2k^{p/2-1}m_{k,p}).$$

Note that (a) and (b) are not contradictory, since $\Omega_t^j = \varnothing$, when $X$ is continuous. It is also worth noticing that the conditional variances (18) and (19), although of course random, are more or less behaving in time like $1/t$.

4.3. *Consistent estimators of the asymptotic variances.* To evaluate the level of tests based on the statistic $\widehat{S}(p,k,\Delta_n)_t$, we need consistent estimators for the asymptotic variances obtained in Theorem 3. That is, we will need to estimate $D(p)$ when $p \geq 2$ and when there are jumps, and also for $A(p)$ when $p \geq 2$ and $X$ is continuous.

To estimate $A(p)$, we can use a realized truncated $p$th variation: for any constants $\alpha > 0$ and $\varpi \in (0, \frac{1}{2})$, we have from Jacod (2008) that, if either $p = 2$, or $p > 2$ and $X$ is continuous, then

$$(21) \qquad \widehat{A}(p,\Delta_n)_t := \frac{\Delta_n^{1-p/2}}{m_p} \sum_{i=1}^{[t/\Delta_n]} |\Delta_i^n X|^p 1_{\{|\Delta_i^n X| \leq \alpha \Delta_n^\varpi\}} \xrightarrow{\mathbb{P}} A(p)_t.$$

Alternatively, we can use the multipower variations of Barndorff-Nielsen and Shephard (2006). For any $r \in (0, \infty)$ and any integer $q \geq 1$, we have from Barndorff-Nielsen et al. (2006a) that if $X$ is continuous

$$(22) \qquad \tilde{A}(r,q,\Delta_n)_t := \frac{\Delta_n^{1-qr/2}}{m_r^q} \sum_{i=1}^{[t/\Delta_n]-q+1} \prod_{j=1}^{q} |\Delta_{i+j-1}^n X|^r \xrightarrow{\mathbb{P}} A(qr)_t$$

[when $q = 1$, we have $\tilde{A}(r,1,\Delta_n) = (\Delta_n^{1-r/2}/m_r)\widehat{B}(r,\Delta_n)$].

It turns out that we also need some results concerning the behavior of $\widehat{A}(p,\Delta_n)_t$ or $\tilde{A}(r,q,\Delta_n)$ when $X$ is discontinuous. These results will be stated below, together with the consistency of the estimators of $D(p)_t$ which we presently describe. Estimators for $D(p)_t$ are a bit more difficult to construct because we need to evaluate $\sigma_{s-}^2$ and $\sigma_s^2$ when $s$ is a jump time, and this involves a kind of nonparametric estimation. A possibility, among many others, is as follows: take any sequence $k_n$ of integers satisfying

$$(23) \qquad\qquad k_n \to \infty, \qquad k_n \Delta_n \to 0$$

and then let $I_{n,t}(i) = \{j \in \mathbb{N} : j \neq i : 1 \leq j \leq [t/\Delta_n], |i-j| \leq k_n\}$ define a local window in time of length $2k_n\Delta_n$ around time $i\Delta_n$ and

$$(24) \quad \widehat{D}(p,\Delta_n)_t = \frac{1}{k_n\Delta_n} \sum_{i=1}^{[t/\Delta_n]} |\Delta_i^n X|^p \sum_{j \in I_{n,t}(i)} (\Delta_j^n X)^2 \mathbf{1}_{\{|\Delta_j^n X| \leq \alpha\Delta_n^\varpi\}},$$

where $\alpha > 0$ and $\varpi \in (0, 1/2)$.

The following theorem establishes the consistency of these estimators, and also contains the technical result (25) whose proof goes along the same lines, and which is needed for the consistency of the tests given below when the null hypothesis is that $X$ does not jump.

THEOREM 4. (a) *We have*

$$(25) \quad \begin{array}{l} p \geq 2, \ t > 0, \ \alpha > 0, \\[6pt] \dfrac{1}{2} - \dfrac{1}{p} < \varpi < \dfrac{1}{2} \quad \Rightarrow \quad \limsup_n \dfrac{\Delta_n \widehat{A}(2p, \Delta_n)_t}{\widehat{A}(p, \Delta_n)_t^2} \xrightarrow{\mathbb{P}} 0, \end{array}$$

$$(26) \quad r \in (0,2), \ q \in \mathbb{N} \quad \Rightarrow \quad \tilde{A}(r,q,\Delta_n)_t \xrightarrow{\mathbb{P}} A(qr)_t,$$

*locally uniformly in $t$.*

(b) *If $\alpha > 0$, $\varpi \in (0, 1/2)$ and $p > 2$, we have*

$$(27) \quad \widehat{D}(p,\Delta_n)_t \xrightarrow{\mathbb{P}} D(p)_t \qquad \forall t \geq 0.$$

*If further $X$ is continuous, then*

$$(28) \quad \Delta_n^{1-p/2} \widehat{D}(p,\Delta_n)_t \xrightarrow{\mathbb{P}} m_p A(p+2)_t \qquad \text{locally uniformly in } t.$$

Note that (27) and (28) are not contradictory, because $D(p) = 0$ when $X$ is continuous. In fact we have more than (25), namely the convergence (21) when there are jumps, but only under some restrictions on the jumps and also some restrictions on $\varpi$ connected with the structure of the jumps and the value $p > 2$; but these refinements are useless here. Of course (26) reduces to (22) when $X$ is continuous.

4.4. *CLT for the standardized statistics.* Using Theorem 4, we can immediately construct consistent estimators of the asymptotic variances of the test statistics established in (18) and (19), respectively. We deduce from (11), (26) and (27) the following CLT for the standardized test statistics [recall $\int_0^t |\sigma_s|^p \, ds > 0$ a.s., by (f) of Assumption 1]:

THEOREM 5. (a) *Let $p > 3$ and $t > 0$. With*

$$\widehat{V}_{n,t}^{j} = \frac{\Delta_n(k-1)p^2 \widehat{D}(2p-2, \Delta_n)_t}{2\widehat{B}(p, \Delta_n)_t^2}, \tag{29}$$

*the variables $(\widehat{V}_{n,t}^{j})^{-1/2}(\widehat{S}(p, k, \Delta_n)_t - 1)$ converge stably in law, in restriction to the set $\Omega_t^j$ of (6), to a variable which, conditionally on $\mathcal{F}$, is centered with variance 1, and which is $N(0, 1)$ if in addition the processes $\sigma$ and $X$ have no common jumps.*

(b) *Assume in addition that $X$ is continuous, and let $p \geq 2$ and $t > 0$. The variables $(\widehat{V}_{n,t}^{c})^{-1/2}(\widehat{S}(p, k, \Delta_n)_t - k^{p/2-1})$ converge stably in law to a variable which, conditionally on $\mathcal{F}$, is $N(0, 1)$, where $\widehat{V}_{n,t}^{c}$ is based on truncations:*

$$\widehat{V}_{n,t}^{c} = \frac{\Delta_n M(p, k) \widehat{A}(2p, \Delta_n)_t}{\widehat{A}(p, \Delta_n)_t^2}, \tag{30}$$

*or is replaced by the multipower estimator:*

$$\tilde{V}_{n,t}^{c} = \frac{\Delta_n M(p, k) \tilde{A}(p/([p]+1), 2[p]+2, \Delta_n)_t}{\tilde{A}(p/([p]+1), [p]+1, \Delta_n)_t^2}. \tag{31}$$

In (31) we have chosen $r = p/([p]+1)$ and respectively $q = 2[p] + 2$ and $q = [p] + 1$. Any other choice with $r \in (0, 2)$ and respectively $q = 2p/r$ and $q = p/r$ would do as well.

**5. Testing for jumps.** We now use the preceding results to construct actual tests, either for the null hypothesis that there are no jumps, or for the null hypothesis that jumps are present.

5.1. *When there are no jumps under the null hypothesis.* In a first case, we set the null hypothesis to be "no jump." We choose an integer $k \geq 2$ and a real $p > 3$, and associate the critical (rejection) region of the form

$$C_{n,t}^{c} = \{\widehat{S}(p, k, \Delta_n)_t < c_{n,t}^{c}\} \tag{32}$$

for some sequence $c_{n,t}^{c} > 0$, possibly $c_{n,t}^{c} = c_t^{c}$ for all $n$, and possibly even a random sequence.

The customary way for defining the asymptotic level of this test is as follows: if $\alpha_{n,t}^{c}(b, \sigma, \delta) = \mathbb{P}(C_{n,t}^{c})$, a notation which emphasizes the dependency on the coefficients $(b, \sigma, \delta)$ of (1), one should take the supremum over all triples $(b, \sigma, \delta)$ in the null hypothesis (i.e., with $\delta \equiv 0$) of $\limsup_n \alpha_t^n(b, \sigma, \delta)$. However, as mentioned at the end of Section 3.1, we only observe a particular path $s \mapsto X_s(\omega)$ over $[0, t]$, and even only at times $i\Delta_n$, so there is obviously no way of statistically separating the genuine null hypothesis from the case where there are jumps, but none occurred in the interval $[0, t]$. Therefore, recalling the sets $\Omega_t^c$ and $\Omega_t^j$ of (6), it is natural to take the following as our definition of the asymptotic level:

$$\alpha = \sup_{b, \sigma, \delta} \limsup_n \mathbb{P}(C_{n,t}^c \mid \Omega_t^c), \tag{33}$$

with the convention that $\mathbb{P}(\cdot \mid \Omega_t^c) = 0$ if $\mathbb{P}(\Omega_t^c) = 0$. Observe that we take first the lim sup and next the supremum; should we proceed the other way around, we would find $\alpha = 1$. In a similar way, the power function for the coefficients triple $(b, c, \delta)$ is the following conditional probability:

$$\beta_{n,t}^c(b, \sigma, \delta) = \mathbb{P}(C_{n,t}^c \mid \Omega_t^j). \tag{34}$$

The right-hand side above makes sense as soon as the function $\delta$ in restriction to $\Omega \times [0,t] \times E$ is not $\mathbb{P}(d\omega)\,ds\,\lambda(dx)$-almost everywhere vanishing, since $\mathbb{P}(\Omega_t^j) > 0$ in that case, whereas otherwise $\mathbb{P}(\Omega_t^j > 0) = 0$ and our process is continuous on $[0, t]$ (since our setting is nonhomogeneous in time, a test based on observations on $[0, t]$ can of course say nothing about what happens after time $t$).

For $\alpha \in (0, 1)$, denote by $z_\alpha$ the $\alpha$-quantile of $N(0,1)$, that is, $\mathbb{P}(U > z_\alpha) = \alpha$, where $U$ is $N(0,1)$. We have:

THEOREM 6. *Let $t > 0$, choose a real $p > 3$ and an integer $k \geq 2$, and let*

$$c_{n,t}^c = k^{p/2-1} - z_\alpha \sqrt{\widehat{V}_{n,t}^c}, \tag{35}$$

*where $\widehat{V}_{n,t}^c$ is given by (30) with $\alpha > 0$ and $1/2 - 1/p < \varpi < 1/2$, or is replaced by $\tilde{V}_{n,t}^c$ of (31). Then:*

*(a) The asymptotic level (33) of the critical region defined by (32) for testing the null hypothesis "no jump" equals $\alpha$.*

*(b) The power function (34) satisfies $\beta_{n,t}^c(b, \sigma, \delta) \to 1$ for all coefficients $(b, \sigma, \delta)$ such that $\mathbb{P}(\Omega_t^j) > 0$ (and of course satisfying Assumption 1).*

5.2. *When jumps are present under the null hypothesis.* In a second case, we set the null hypothesis to be that "there are jumps," in which case the critical (rejection) region is of the form

$$C_{n,t}^j = \{\widehat{S}(p, k, \Delta_n)_t > c_{n,t}^j\} \tag{36}$$

for some sequence $c_{n,t}^j > 0$. As in (33), the asymptotic level is

$$\alpha' = \sup_{b, \sigma, \delta} \limsup_n \mathbb{P}(C_{n,t}^j \mid \Omega_t^j), \tag{37}$$

with the convention that $\mathbb{P}(\cdot \mid \Omega_t^j) = 0$ if $\mathbb{P}(\Omega_t^j) = 0$, and the power function for the coefficients triple $(b, \sigma, \delta)$ is

$$\beta_{n,t}^j(b, \sigma, \delta) = \mathbb{P}(C_{n,t}^j \mid \Omega_t^c). \tag{38}$$

The right-hand side above is simply $\mathbb{P}(C_{n,t}^j)$ when $\delta = 0$, and it makes sense as soon as $\mathbb{P}(\Omega_t^c) > 0$, whereas when $\mathbb{P}(\Omega_t^c) = 0$, then the null hypothesis is of course satisfied.

THEOREM 7. *Let $t > 0$ and choose a real $p > 3$ and an integer $k \geq 2$. Let $\alpha \in (0, 1)$ and $\varpi \in (0, 1/2)$.*

(a) *Letting*

$$c_{n,t}^j = 1 + \sqrt{\widehat{V}_{n,t}^j / \alpha}, \tag{39}$$

*where $\widehat{V}_{n,t}^j$ is given by (29), the asymptotic level (37) of the critical region defined by (36) for testing the null hypothesis "there are jumps" is smaller than $\alpha$.*

(b) *Suppose that we restrict our attention to models in which the processes $X$ and $\sigma$ have no common jumps. If*

$$c_{n,t}^j = 1 + z_\alpha \sqrt{\widehat{V}_{n,t}^j}, \tag{40}$$

*then the asymptotic level (37) of the critical region defined by (36) for testing the null hypothesis "there are jumps" is equal to $\alpha$.*

(c) *In all cases the power function (38) satisfies $\beta_{n,t}^j(b, \sigma, \delta) \to 1$ for all coefficients $(b, \sigma, \delta)$ such that $\mathbb{P}(\Omega_t^c) > 0$.*

5.3. *Choosing the free parameters $p$, $k$, $\alpha$ and $\varpi$.* In both Theorems 6 and 7 one has two "basic" parameters to choose, namely $p > 3$ and the integer $k \geq 2$. For the standardization, and except if we use $\widetilde{V}_{n,t}^c$ in Theorem 6, we also need to choose $\alpha > 0$ and $\varpi$, which should be in $(1/2 - 1/p, 1/2)$ for the first theorem, and in $(0, 1/2)$ for the second one. Let us give some remarks regarding the selection of these parameters, and the impact of the choices on the test.

First, about the real $p > 3$: the larger it is, the more the emphasis put on "large" jumps. Since those are in any case relatively easy to detect, or at least much easier than the small ones, it is probably wise to choose $p$ "small." However, choosing $p$ close to 3 may induce a rather poor fit in the central limit theorem for a number $n$ of observations which is not very large, since for $p = 3$ there is still a CLT, but with a bias. A good compromise seems to be $p = 4$, since in addition the computations are easier when $p$ is an integer.

Second, about $k$: when $k$ increases, we have to separate two points (1 and $k$ when $p = 4$) which are further and further apart, but on the other hand in

(18) and (19) we see that the asymptotic variances are increasing with $k$. Furthermore, large values of $k$ lead to a decrease in the effective sample size employed to estimate the numerator $\widehat{B}(p, k\Delta_n)_t$ of the test statistic, which is inefficient. So one should choose $k$ not too big. Numerical experiments with $k = 2, 3, 4$ have been conducted below, showing no significant differences for the results of the test itself. More experiments should probably be conducted; however, we think that the choice $k = 2$ is in all cases reasonable.

Third, about $\alpha$ and $\varpi$: from a number of numerical experiments, not reported here to save space, it seems that choosing $\varpi$ close to $1/2$ (as $\varpi = 0.47$ or $\varpi = 0.48$), and $\alpha$ between 3 and 5 times the "average" value of $\sigma$, leads to the left-hand sides of (27) and (28) being very close to the right-hand sides, for relatively small values of $n$. Of course choosing $\alpha$ as above may seem circular because $\sigma_t$ is unknown, and usually random and time-varying, but in practice one very often has a pretty good idea of the order of magnitude of $\sigma_t$, especially for financial data. Even if no a priori order of magnitude is available, it is possible to estimate consistently the volatility of the continuous part of the semimartingale, $(\int_0^t \sigma_s^2 \, ds)^{1/2}$, in the presence of jumps; see the literature on disentangling jumps from diffusions cited in the Introduction. The multipower variations (22) do not suffer from the drawback of having to choose $\alpha$ and $\varpi$ a priori, but they cannot be used for Theorem 7, and when there are jumps the quality of the approximation in (26) strongly depends on the relative sizes of $\sigma$ and of the cumulated jumps.

Finally, let us calculate more explicitly the critical regions for the values $p = 4$. We obtain that $M(4, 2) = 160/3$, and more generally $M(4, k) = 16k(2k^2 - k - 1)/3$, and the cut-off point in (35) when further $k = 2$ becomes

$$
\begin{aligned}
c_{n,t}^c &= 2 - z_\alpha \sqrt{\frac{160\Delta_n \widehat{A}(8, \Delta_n)_t}{3\widehat{A}(4, \Delta_n)_t^2}} \quad \text{or} \\
c_{n,t}^c &= 2 - z_\alpha \sqrt{\frac{160\Delta_n \tilde{A}(4/5, 10, \Delta_n)_t}{3\tilde{A}(4/5, 5, \Delta_n)_t^2}}.
\end{aligned}
\tag{41}
$$

In a similar way, (40) becomes

$$
c_{n,t}^j = 1 + z_\alpha \sqrt{\frac{8\Delta_n \widehat{D}(6, \Delta_n)_t}{\widehat{B}(4, \Delta_n)_t^2}}.
\tag{42}
$$

The $\alpha$-quantiles of $N(0, 1)$ at the 10% and 5% level are $z_{0.1} = 1.28$ and $z_{0.05} = 1.64$, respectively.

5.4. *The effect of microstructure noise.* As already said before, the observations of the process $X$ are blurred with a (small) noise, which messes things up when data are recorded with high-frequency.

We do not intend to provide a deep study of this topic here, but merely to establish some basic facts from the point of view of the consistency of our statistics (rates of convergence are more difficult to obtain). We assume that each observation is affected by an additive noise, that is, instead of $X_{i\Delta_n}$ we really observe $Y_{i\Delta_n} = X_{i\Delta_n} + \varepsilon_i$, and the $\varepsilon_i$ are supposed to be i.i.d. with $E(\varepsilon_i^2)$ and $E(\varepsilon_i^4)$ finite. Then, instead of $\widehat{B}(4, \Delta_n)_t$ (we consider the case $p = 4$ here), we actually observe (for $k = 1$ and also $k \geq 2$):

$$\widehat{B}'(4, k\Delta_n)_t = \sum_{i=1}^{[t/k\Delta_n]} (X_{ik\Delta_n} - X_{(i-1)k\Delta_n} + \varepsilon_{ki} - \varepsilon_{k(i-1)})^4$$

$$= \widehat{B}(4, k\Delta_n)_t + 2 \sum_{i=1}^{[t/k\Delta_n]} (X_{ik\Delta_n} - X_{(i-1)k\Delta_n})^2 (\varepsilon_{ki} - \varepsilon_{k(i-1)})^2$$

$$+ \sum_{i=1}^{[t/k\Delta_n]} (\varepsilon_{ki} - \varepsilon_{k(i-1)})^4.$$

The second term above behaves like $4E(\varepsilon_i^2)\widehat{B}(2, k\Delta_n)_t$ and the third one like $t(2E(\varepsilon_i^4) + 6E(\varepsilon_i^2)^2)/(k\Delta_n)$. It follows that the statistics $\widehat{S}'(4, k, \Delta_n)_t = \widehat{B}'(4, k\Delta_n)_t / \widehat{B}'(4, \Delta_n)_t$ that we actually observe instead of (12) has the following behavior, as $\Delta_n \to 0$:

$$(43) \qquad \widehat{S}'(4, k, \Delta_n)_t \xrightarrow{\mathbb{P}} \frac{1}{k}.$$

The relevance of this limit will become clear when we apply the test to real data in Section 7 below. Finally, note that when $\Delta_n$ is moderately small, things may be different: if $E(\varepsilon_i^2)$ and $E(\varepsilon_i^4)$ are small (as they are in practice), $\widehat{S}'(4, k, \Delta_n)_t$ will be close to 1 on $\Omega_t^j$ and to $k$ on $\Omega_t^c$.

**6. Simulation results.** Throughout the simulations, we use $p = 4$. We calibrate the values to be realistic for a liquid stock trading on the NYSE. We use an observation length of $t = 1$ day, consisting of 6.5 hours of trading, that is, 23,400 seconds. The simulations contain no microstructure noise.

Under the null of no jumps, the performance of the test in simulations is reported in Table 1. The distribution of the test statistic in simulations is close to the theoretical normal limit centered at $k$ in simulations; histograms for the test statistic and the corresponding asymptotic distribution are shown in Figure 1 for the nonstandardized (top panel) and standardized test statistics (lower panel) and the cases $k = 2$ (left panel) and $k = 3$

TABLE 1
*Level of the test under the null hypothesis of no jumps*

| $\Delta_n$ | $n$ | $k$ | Mean value of $\widehat{S}(4, k, \Delta_n)$ | | Rejection rate in simulations | |
|---|---|---|---|---|---|---|
| | | | Asymptotic | Simulations | 10% | 5% |
| 1 sec | 23,400 | 2 | 2 | 2.000 | 0.098 | 0.046 |
| 1 sec | 23,400 | 3 | 3 | 2.998 | 0.099 | 0.048 |
| 1 sec | 23,400 | 4 | 4 | 3.999 | 0.094 | 0.043 |
| 5 sec | 4680 | 2 | 2 | 1.998 | 0.099 | 0.044 |
| 5 sec | 4680 | 3 | 3 | 2.998 | 0.099 | 0.043 |
| 5 sec | 4680 | 4 | 4 | 3.995 | 0.095 | 0.038 |
| 10 sec | 2340 | 2 | 2 | 2.000 | 0.093 | 0.041 |
| 15 sec | 1560 | 2 | 2 | 1.999 | 0.090 | 0.043 |
| 30 sec | 780 | 2 | 2 | 2.007 | 0.085 | 0.035 |

Note: This table reports the results of 5000 simulations of the test statistic under the null hypothesis of no jumps. The data generating process is the stochastic volatility model $dX_t/X_t = \sigma_t \, dW_t$, with $\sigma_t = v_t^{1/2}$, $dv_t = \kappa(\beta - v_t)\,dt + \gamma v_t^{1/2}\,dB_t$, $E[dW_t\,dB_t] = \rho\,dt$, $\beta^{1/2} = 0.4$, $\gamma = 0.5$, $\kappa = 5$, $\rho = -0.5$. The parameter values are realistic for a stock based on Aït-Sahalia and Kimmel (2007). The test statistic is standardized with the estimator of $\widehat{V}_{n,t}^c$ given in (30) with $\alpha = 5\beta^{1/2}$ and $\varpi = 0.47$ in (21).

(right panel.) The differences between the cases $k = 2$ and $k = 3$ conform to the theoretical arguments given above in Section 5.3, in terms of trade-off between variance and separation of the modes of the distributions.

Under the null that jumps are present, the test statistic is now centered around the predicted value of 1. We report in Table 2 the performance of the test statistic when jumps are compound Poisson. Figure 2 plots the distribution of the test statistic for different values of the jump arrival rate, $\lambda$. An interesting phenomenon happens, when multiple jumps can occur with high probability ($\lambda$ high). While the likelihood that jumps will take place in successive observations is small, such paths will happen over a large number of simulations. When that situation happens, we may see in $\widehat{B}(4, 2\Delta_n)_t$ the two jumps either compensate each other, if they are of opposite signs, or cumulate into a single larger jump if they are of the same sign. On the other hand, no such compensation or cumulation takes place in $\widehat{B}(4, \Delta_n)_t$. As a result, their ratio $\widehat{S}(4, 2, \Delta_n)_t$ can exhibit a small number of outliers.

Similar results for the case where jumps are generated by a Cauchy process are reported in Table 3 and in Figure 3 for the nonstandardized (top panel) and standardized test statistics (bottom panel), for different values of the Cauchy scale parameter $\theta$. The contrast between the top (very far from normal) and bottom (nearly normal) panels in Figure 3 illustrates the role of the standardization in making the test statistic asymptotically normal: recall from Theorem 3(a) that the nonstandardized statistic $\widehat{S}(4, 2, \Delta_n)_t$ is in general non-Normal unconditionally, whereas from Theorem 5(a) the standardized statistic $(\widehat{V}_{n,t}^j)^{-1/2}(\widehat{S}(4,2,\Delta_n)_t - 1)$ is asymptotically $N(0,1)$ unconditionally.

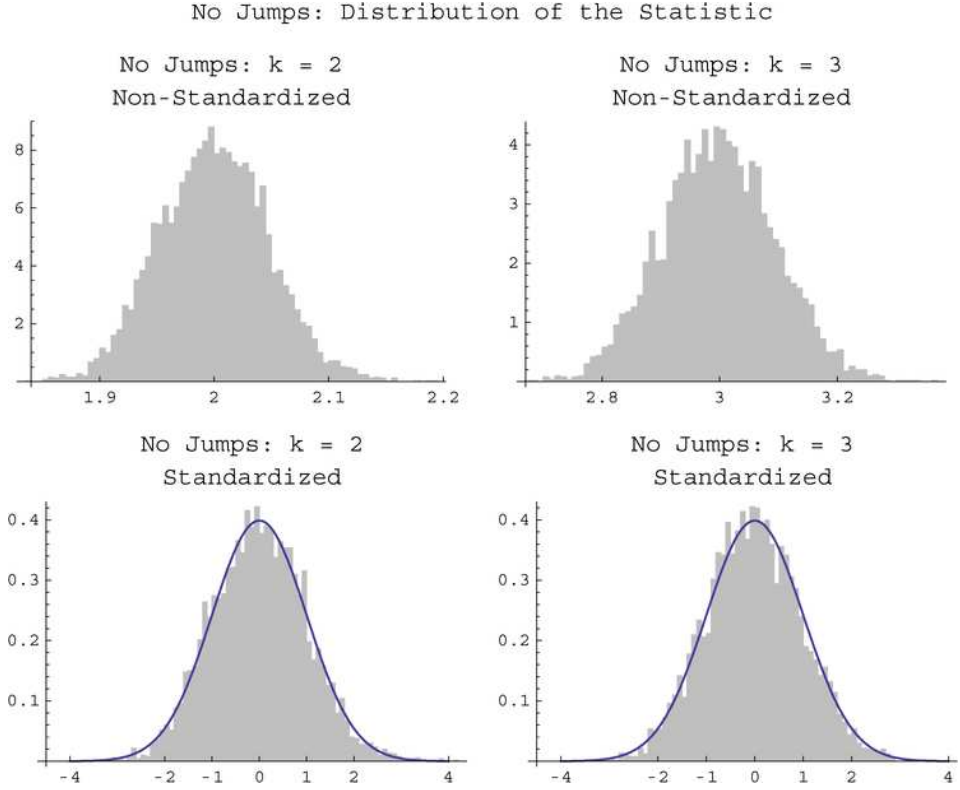

Fig. 1. *Monte Carlo and theoretical asymptotic distributions of the non-standardized (top row) and standardized (bottom row) test statistics $S(4,k,\Delta_n)$ for $k=2,3$ and $\Delta_n = 1$ second under the null hypothesis of no jumps. In the standardized case, the solid curve is the $N(0,1)$ density.*

In Figure 4, we illustrate what happens when there are paths containing either no jumps or jumps that are too tiny to be detected as jumps. In the Poisson case (left panel), we set $\lambda$ to correspond to 1 jump per day on average, but keep all paths, including those where no jump took place. Such paths with no jumps occur with positive probability, unlike the Cauchy case. But since we condition on having jumps, we removed those paths in Figure 2 to compute the distribution of the statistic under the null that jumps are present, as required under Theorem 5(a); this is the practical implication of the restriction to the set $\Omega_t^j$. Now, when those paths are kept, we obtain a clear-cut bimodal result: either no jump occurred, and those

TABLE 2
*Compound Poisson jumps: Level of the test under the null hypothesis that jumps are present*

| | | | Mean value of $\widehat{S}(4,k,\Delta_n)$ | | Rejection rate in simulations | |
|---|---|---|---|---|---|---|
| $\Delta_n$ | $n$ | $k$ | Asymptotic | Simulations | 10% | 5% |
| | | | $\lambda = 1$ jump per day | | | |
| 1 sec | 23,400 | 2 | 1 | 1.000 | 0.110 | 0.056 |
| 4 sec | 5,850 | 2 | 1 | 1.002 | 0.107 | 0.054 |
| 15 sec | 1560 | 2 | 1 | 1.003 | 0.100 | 0.053 |
| | | | $\lambda = 5$ jumps per day | | | |
| 1 sec | 23,400 | 2 | 1 | 1.002 | 0.113 | 0.057 |
| 4 sec | 5850 | 2 | 1 | 1.006 | 0.112 | 0.059 |
| 15 sec | 1560 | 2 | 1 | 1.012 | 0.129 | 0.071 |
| | | | $\lambda = 10$ jumps per day | | | |
| 1 sec | 23,400 | 2 | 1 | 1.002 | 0.105 | 0.056 |
| 4 sec | 5850 | 2 | 1 | 1.009 | 0.134 | 0.076 |
| 15 sec | 1560 | 2 | 1 | 1.029 | 0.163 | 0.099 |

Note: This table reports the results of 5000 simulations of the test statistic under the null hypothesis that jumps are present. The model under the null is $dX_t/X_t = \sigma_t\,dW_t + J_t\,dN_t$, where $\sigma_t$ is the same stochastic volatility process with the same parameter values as in Table 1, $J_t$ is the product of a uniformly distributed variable on $[-2,-1] \cup [1,2]$ times a constant $J_S$ and $N$ is a Poisson process with intensity $\lambda$. The total variance of the increments, $\sigma^2 + (7/3)J_S^2\lambda$ is held constant at $0.4^2$. As a result, jumps that are more frequent tend to be smaller in size. In the simulations, 25% of the total variance is due to the Brownian motion and 75% to the jumps. Since the test is conditional on a path containing jumps, paths that do not contain any jump are excluded from the simulated sample and replaced by new simulations. Thus, in sample, the number of jumps is slightly higher than specified by the value of $\lambda$ in the table. The test statistic is standardized with the estimator of $\widehat{V}_{n,t}^j$ given in (29) with $k_n = [50\Delta_n^{-1/4}]$ in (24).

paths are grouped around 2, or (at least) 1 jump occurred and those paths are grouped around 1.

In the Cauchy case (right panel of Figure 4), when a small value of $\theta$ is selected, although each path has infinitely many jumps, a sizeable proportion of paths has no jump large enough to make the path look markedly different from that of a pure Brownian motion at our observation frequency. As a result, we get a bimodal distribution with a second mode at the (continuous case) value of 2. Compared to the left panel of the figure, the infinite-activity of the process produces a more diffuse situation where some paths have only tiny jumps, and the statistic takes intermediary values between 1 and 2.

Finally, we confirm in simulations that jumps in $\sigma$ do not affect the distribution of the test statistic; as predicted by Theorem 5, this is always the case if $X$ is continuous, and when $X$ jumps, remains the case as long as $\sigma$

Poisson Jumps: Distribution of the Statistic

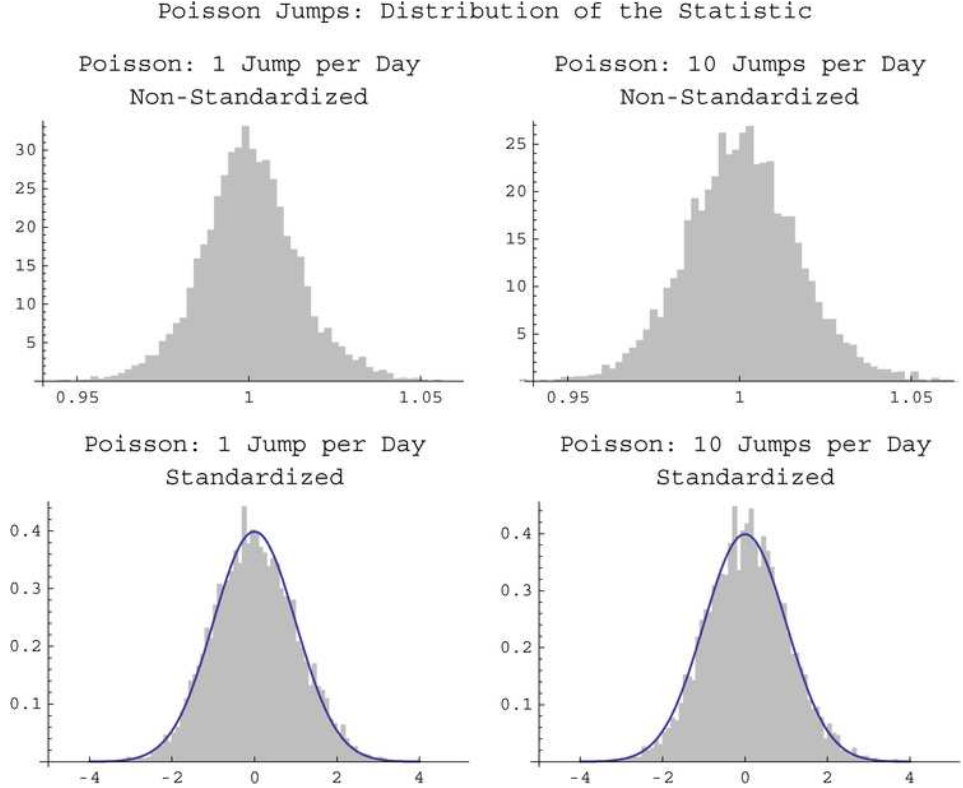

FIG. 2. *Monte Carlo and theoretical $N(0,1)$ asymptotic distributions of the standardized test statistic $(\widehat{V}_n^j)^{-1/2}(\widehat{S}(4,2,\Delta_n)-1)$ for $\Delta_n = 1$ second under the null hypothesis that jumps are present, using the same data generating process with Poisson jumps as in Table 2.*

and $X$ have no common jumps. To check this, we repeat the simulations above with jumps added to $\sigma$ that are independent from those of $X$ (if any): we simulate $v = \sigma^2$ from the model used in Table 1 plus proportional compound Poisson jumps that are uniformly distributed on $[-30\%, 30\%]$. The results are largely unchanged from those in the tables above, and are not reported here to save space.

**7. Empirical application.** We now conduct the test for each of the 30 Dow Jones Industrial Average (DJIA) stocks and each trading day in 2005; the data source is the TAQ database. Each day, we collect all transactions on the NYSE or NASDAQ, from 9:30am until 4:00pm, for each one of these stocks. We sample in calendar time every 5 seconds. Each day and stock is treated on its own. We use filters to eliminate clear data errors (price set to zero, etc.) as is standard in the empirical market microstructure literature.

TABLE 3
*Cauchy jumps: level of the test under the null hypothesis that jumps are present*

| $\Delta_n$ | $n$ | $k$ | Mean value of $\widehat{S}(4,k,\Delta_n)$ | | Rejection rate in simulations | |
|---|---|---|---|---|---|---|
| | | | Asymptotic | Simulations | 10% | 5% |
| | | | | Jump size: $\theta = 10$ | | |
| 1 sec | 23,400 | 2 | 1 | 1.002 | 0.112 | 0.062 |
| 5 sec | 4680 | 2 | 1 | 1.010 | 0.117 | 0.066 |
| 15 sec | 1560 | 2 | 1 | 1.027 | 0.126 | 0.081 |
| | | | | Jump size: $\theta = 50$ | | |
| 1 sec | 23,400 | 2 | 1 | 1.000 | 0.090 | 0.051 |
| 5 sec | 4680 | 2 | 1 | 1.002 | 0.100 | 0.059 |
| 15 sec | 1560 | 2 | 1 | 1.003 | 0.094 | 0.059 |

Note: This table reports the results of 5000 simulations of the test statistic under the null hypothesis that jumps are present. The model under the null is $dX_t/X_t = \sigma_t\,dW_t + \theta\,dY_t$, where $\sigma_t$ is the same stochastic volatility process as in Tables 1 and 2, $Y$ is a Cauchy process standardized to have characteristic function $\mathbb{E}(\exp(iuY_t)) = \exp(-t|u|/2)$. The value of the long-run mean of volatility, $\beta^{1/2}$, is 0.2; the other parameter values are identical. Given $\beta$, the parameter $\theta$ measures the size of the jumps relative to the volatility. The test statistic is standardized with the estimator of $\widehat{V}_{n,t}^j$ given in (29) with $k_n = [50\Delta_n^{-1/4}]$ in (24).

We plot in Figure 5 a histogram showing the empirical distribution of the statistic computed on these data. As expected, we see evidence of market microstructure noise in the form of density mass below 1 (the limit is 1/2, see Section 5.4, if the noise is i.i.d., but may be different with other kinds of noise). The striking feature of the results is that most of the observed values are around 1, providing evidence for the presence of jumps, with only few observations around 2, the expected limit in the continuous case. As we sample less frequently, the distribution spreads out, consistently with both the asymptotic theory and the simulations above.

**8. Proofs.**

8.1. *Preliminaries.* We use the shorthand notation $E_{i-1}^n(Y)$ for $E(Y \mid \mathcal{F}_{(i-1)\Delta_n})$, and we set

$$\delta_i^n = \sigma_{(i-1)\Delta_n}\Delta_i^n W, \qquad \theta_i^n = |\Delta_i^n X|1_{\{|\Delta_i^n X|\leq \alpha\Delta_n^\varpi\}}$$

for a given pair $\alpha > 0$ and $\varpi \in (0, \frac{1}{2})$. We will consider a strengthened version of Assumption 1:

ASSUMPTION 2. We have Assumption 1, and $|b_t| + |\sigma_t| + |\widetilde{b}_t| + |\widetilde{\sigma}_t| + |\widetilde{\sigma}_t'| \leq K$ and $|\delta(t,x)| \leq \gamma(x)$ and $|\widetilde{\delta}(t,x)| \leq \gamma(x)$ and also $\gamma(x) \leq K$ for some constant $K$.

TESTING FOR JUMPS 23

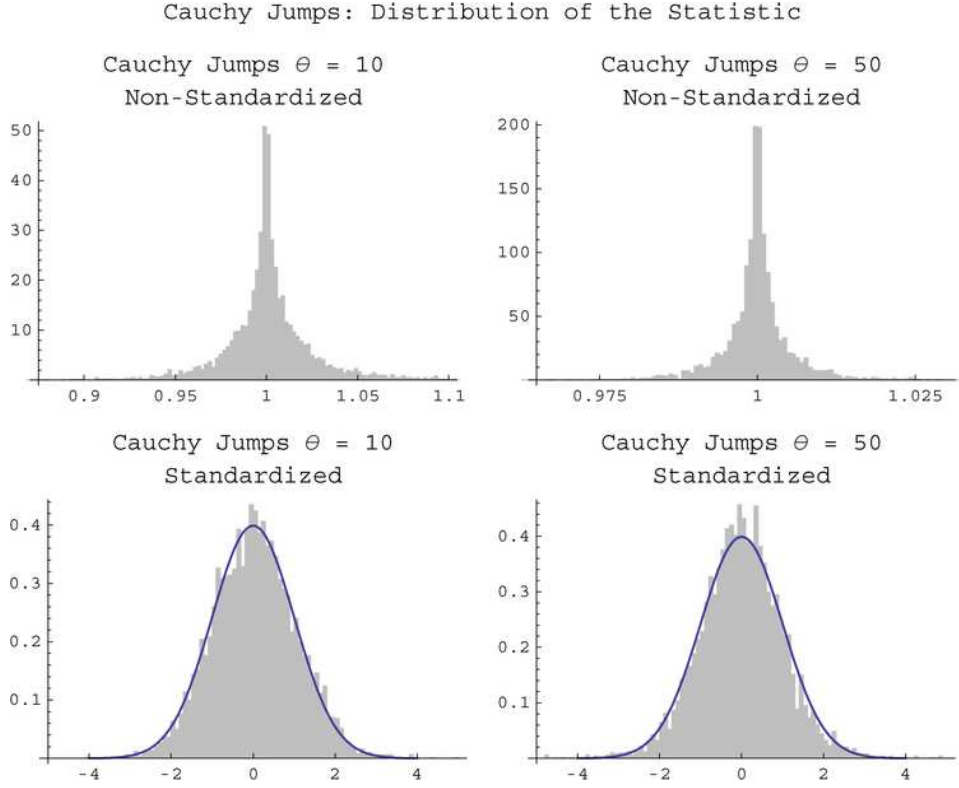

Cauchy Jumps: Distribution of the Statistic

FIG. 3. *Monte Carlo distribution of the nons-tandardized test statistic $S(4,2,\Delta_n)$ for $\Delta_n = 1$ second under the null hypothesis that jumps are present, using the same data generating process with Cauchy jumps as in Table 3. The asymptotic distribution of the nonstandardized statistic is nonnormal as expected (see upper panel). However, the asymptotic distribution of the standardized test statistic is $N(0,1)$ (the solid curve in the lower panel).*

A localization procedure shows that for proving Theorems 2 and 4 we can replace everywhere Assumption 1 by Assumption 2; this procedure is described in detail in Jacod (2008) and works with no change at all here, so we omit it.

Now we state a number of consequences of this strengthened assumption, to be used a number of times in the following proofs. Below, $K$ denotes a constant which may depend on the coefficients $(b,c,\delta)$ and $(\widetilde{b},\widetilde{\sigma},\widetilde{\sigma}',\widetilde{\delta}')$ and which changes from line to line. We write it $K_a$ if it depends on an additional parameter $a$. We also write $X'_t$ and $X''_t$ for the sum of the first three terms, resp. of the last two terms, in (1). First, for all $q \geq 2$ we have the following classical estimates [proved, e.g., in Jacod (2008)] under Assumption 2:

$$\mathbb{E}^n_{i-1}(|\Delta^n_i X'|^q + |\delta^n_i|^q) \leq K_q \Delta^{q/2}_n, \qquad \mathbb{E}^n_{i-1}(|\Delta^n_i X''|^q) \leq K_q \Delta_n,$$

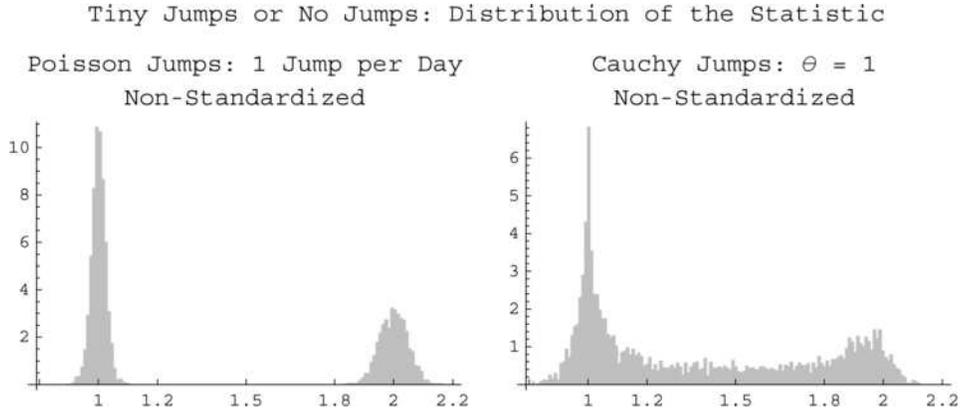

Fig. 4. *Monte Carlo distribution of the non-standardized test statistic $S(4,2,\Delta_n)$ for $\Delta_n = 1$ second, computed using a data generating process with either one Poisson jump per day on average including paths that contain no jumps (left panel) or tiny Cauchy jumps ($\theta = 1$, right panel).*

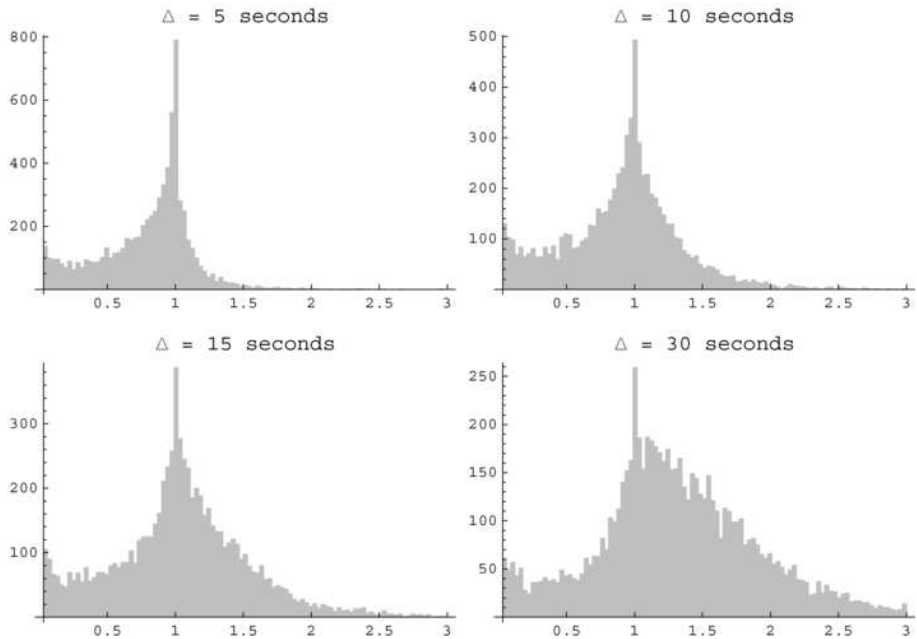

Fig. 5. *Empirical distribution of the nonstandardized statistic $\widehat{S}(4,2,\Delta_n)$ for different values of the sampling interval $\Delta_n$. Each sample point is computed using all the transactions for one of the 30 DJIA stocks observed over one trading day in 2005. This produces 7560 realizations of the statistic.*

(44)
$$\mathbb{E}^n_{i-1}(|\Delta^n_i\sigma|^q) \leq K_q\Delta_n, \qquad \mathbb{E}^n_{i-1}(|\Delta^n_iX' - \delta^n_i|^q) \leq K\Delta_n^{1+q/2}.$$

Next, the proof of Lemma 5.12 of Jacod (2008) applied with $\varepsilon$ being the supremum of the bounded function $\gamma$ shows that under Assumption 2,

(45)
$$\mathbb{E}^n_{i-1}(|\Delta^n_iX''|^2 \wedge \eta^2) \leq K\Delta_n\left(\frac{\eta^2+\Delta_n}{\theta^2} + \Gamma(\theta)\right)$$

for all $\eta > 0$ and $\theta \in (0,1)$, and where $\Gamma(\theta) \to 0$ as $\theta \to 0$. Combining this with the second line in (44) yields that also

(46)
$$\mathbb{E}^n_{i-1}(|\Delta^n_iX - \delta^n_i|^2 \wedge \eta^2) \leq K\Delta_n\left(\frac{\eta^2+\Delta_n}{\theta^2} + \Gamma(\theta)\right).$$

Finally let us also mention two inequalities, valid for all $x,y \in \mathbb{R}$ and $0 < \varepsilon < 1 < A$ and for $p \geq 2$ and $0 < r < 2$:

(47)
$$\begin{aligned}&||x+y|^p 1_{\{|x+y|\leq\varepsilon\}} - |x|^p| \\ &\quad \leq K_p(|x|^p 1_{\{|x|>\varepsilon/2\}} + \varepsilon^{p-2}(y^2 \wedge \varepsilon^2) + |x|^{p-1}(|y| \wedge \varepsilon))\end{aligned}$$

and

(48)
$$\begin{aligned}&||x+y|^r - |x|^r| \\ &\quad \leq K_r(\varepsilon^r + A\varepsilon + A^{r-2}(x^2+y^2) + A^r\varepsilon^{-2}(y^2 \wedge 1)).\end{aligned}$$

These inequalities are elementary, although a bit tedious to show: for the first one, one singles out three cases, namely the case $|x| > \varepsilon/2$, the case $|x| \geq \varepsilon/2$ and $|y| \geq \varepsilon/2$, and the case $|x| \leq \varepsilon/2 < |y|$, and proves the inequality in each of the cases. The second one is proved analogously, after singling out five cases: the case $|x| \leq |y|$, the case $|x| > |y|$ and $|x| > A$, the case $|y| < |x| \leq \varepsilon$, the case $|y| \leq \varepsilon < |x| \leq A$, and the case $\varepsilon < |y| < |x| \leq A$.

8.2. *Proof of Theorem 2.* For proving Theorem 2 we need to exhibit the limits $Z(p)$, $Z'(p,k)$, $Y(p)$ and $Y'(p,k)$ and it takes some additional notation to do so. We consider an auxiliary space $(\Omega', \mathcal{F}', \mathbb{P}')$ which supports a number of variables and processes:

- four sequences $(U_q)$, $(U'_q)$, $(\overline{U}_q)$, $(\overline{U}'_q)$ of $N(0,1)$ variables;
- a sequence $(\kappa_q)$ of uniform variables on $[0,1]$;
- a sequence $(L_q)$ of uniform variables on the finite set $\{0,1,\ldots,k-1\}$ ($k \geq 2$ is the integer showing up in the theorem);
- two standard Wiener processes $\overline{W}$ and $\overline{W}'$;

and all these processes or variables are mutually independent. Then we put

$$\widetilde{\Omega} = \Omega \times \Omega', \qquad \widetilde{\mathcal{F}} = \mathcal{F} \otimes \mathcal{F}', \qquad \widetilde{\mathbb{P}} = \mathbb{P} \otimes \mathbb{P}'$$

and we extend the variables $X_t, b_t, \ldots$ defined on $\Omega$ and $\overline{W}, U_q, \ldots$ defined on $\Omega'$ to the product $\widetilde{\Omega}$ in the obvious way, without changing the notation. We write $\widetilde{\mathbb{E}}$ for the expectation w.r.t. $\widetilde{\mathbb{P}}$. Finally, denote by $(T_n)_{n \geq 1}$ an enumeration of the jump times of $X$ which are stopping times, and let $(\widetilde{\mathcal{F}}_t)$ be the smallest (right-continuous) filtration of $\widetilde{\mathcal{F}}$ containing the filtration $(\mathcal{F}_t)$ and w.r.t. which $\overline{W}$ and $\overline{W}'$ are adapted and such that $U_n, U_n', \overline{U}_n, \overline{U}_n', \kappa_n$ and $L_n$ are $\widetilde{\mathcal{F}}_{T_n}$-measurable for all $n$. We hence get an extension $(\widetilde{\Omega}, \widetilde{\mathcal{F}}, (\widetilde{\mathcal{F}}_t)_{t \geq 0}, \widetilde{\mathbb{P}})$ of the original space $(\Omega, \mathcal{F}, (\mathcal{F}_t)_{t \geq 0}, \mathbb{P})$. Obviously, $W$, $W'$, $\overline{W}$, $\overline{W}'$ are Wiener processes and $\mu$ is a Poisson measure with compensator $\nu$ on $(\widetilde{\Omega}, \widetilde{\mathcal{F}}, (\widetilde{\mathcal{F}}_t)_{t \geq 0}, \widetilde{\mathbb{P}})$.

Now we exhibit the limits. If $(h_s)$ and $(h_s')$ are adapted processes with right-continuous or left-continuous paths, defined on $(\Omega, \mathcal{F}, (\mathcal{F}_t)_{t \geq 0}, \mathbb{P})$, we set

$$Y(h, h')_t = \int_0^t h_s \, d\overline{W}_s + \int_0^t h_s' \, d\overline{W}_s'. \tag{49}$$

This defines a local martingale on the extension, which conditionally on the $\sigma$-field $\mathcal{F}$ is a centered Gaussian martingale. Furthermore if $(k_s, k_s')$ is another pair of processes we have

$$\widetilde{\mathbb{E}}(Y(h, h')_t Y(k, k')_t \mid \mathcal{F}) = \int_0^t (h_s k_s + h_s' k_s') \, ds \tag{50}$$

as in Jacod (2008). In a similar way, with any function $g$ on $\mathbb{R}$ which is locally bounded and with $g(x)/x \to 0$ as $x \to 0$, we associate the following two processes:

$$Z(g)_t = \sum_{q : T_q \leq t} g(\Delta X_{T_q})(\sqrt{\kappa_q} U_q \sigma_{T_q-} + \sqrt{1 - \kappa_q} U_q' \sigma_{T_q}), \tag{51}$$

$$Z'(g)_t = \sum_{q : T_q \leq t} g(\Delta X_{T_q})(\sqrt{L_q} \overline{U}_q \sigma_{T_q-} + \sqrt{k - 1 - L_q} \overline{U}_q' \sigma_{T_q}). \tag{52}$$

That $Z(g)$ is well defined is proved in Jacod (2008) and the proof that $Z'(g)$ is also well defined is exactly the same. Again, conditionally on $\mathcal{F}$, these two processes are independent martingales with mean 0 and

$$\widetilde{\mathbb{E}}(Z(g)_t^2 \mid \mathcal{F}) = \tfrac{1}{2} \sum_{s \leq t} g(\Delta X_s)^2 (\sigma_{s-}^2 + \sigma_s^2),$$
$$\widetilde{\mathbb{E}}(Z'(g)_t^2 \mid \mathcal{F}) = \frac{k-1}{2} \sum_{s \leq t} g(\Delta X_s)^2 (\sigma_{s-}^2 + \sigma_s^2). \tag{53}$$

Moreover, conditionally on $\mathcal{F}$, $Y(h, h')$ and $(Z(g), Z'(g))$ are independent, and the latter is also a Gaussian martingale as soon as the processes $X$ and $\sigma$ have no common jumps.

At this stage we can proceed to proving Theorem 2(a). As a matter of fact, it is clearer to prove a more general result: with $k$ an integer strictly bigger than 1 and with any function $f$ on $\mathbb{R}$, we consider the two processes

(54)
$$V^n(f)_t = \sum_{i=1}^{[t/\Delta_n]} f(\Delta_i^n X),$$

$$\overline{V}^n(f)_t = \sum_{i=1}^{[t/k\Delta_n]} f(\Delta_{ki-k+1}^n X + \Delta_{ki-k+2}^n X + \cdots + \Delta_{ki}^n X).$$

We also write $V(f)_t = \sum_{s \leq t} f(\Delta X_s)$, which is well defined as soon as $f(x)/x^2$ is bounded on a neighborhood of 0.

THEOREM 8. *Under Assumption 1, let $f$ be $C^2$ with $f(0) = f'(0) = f''(0) = 0$ ($f'$ and $f''$ are the first two derivatives). The pair of processes*

$$\left(\frac{1}{\sqrt{\Delta_n}}(V^n(f)_t - V(f)_{\Delta_n[t/\Delta_n]}), \frac{1}{\sqrt{\Delta_n}}(\overline{V}^n(f)_t - V(f)_{k\Delta_n[t/k\Delta_n]})\right)$$

*converges stably in law, on the product $\mathbb{D}(\mathbb{R}_+, \mathbb{R}) \times \mathbb{D}(\mathbb{R}_+, \mathbb{R})$ of the Skorokhod spaces, to the process $(Z(f'), Z(f') + Z'(f'))$.*

We have the (stable) convergence in law of the above processes, as elements of the product functional space $\mathbb{D}(\mathbb{R}_+, \mathbb{R})^2$, but usually not as elements of the space $\mathbb{D}(\mathbb{R}_+, \mathbb{R}^2)$ with the (two-dimensional) Skorokhod topology, because a jump of $X$ entails a jump for both components above, but "with a probability close to $j/k$" the times at which these two components jump differ by an amount $j\Delta_n$, for $j = 1, \ldots, k-1$.

PROOF. As said in the beginning of this section, we can assume Assumption 2. The proof is an extension of the proof of Theorem 2.12(i) of Jacod (2008). We start the proof under the additional assumption that $f$ vanishes on $[-2k\varepsilon, 2k\varepsilon]$ for some $\varepsilon > 0$. Let $S_q$ be the successive jump times of the Poisson process $\mu([0,t] \times \{x : \gamma(x) > \varepsilon\})$. Let $\Omega_n(t, \varepsilon)$ be the set of all $\omega$ such that each interval $[0,t] \cap (i\Delta_n, (i+k)\Delta_n]$ contains at most one $S_q(\omega)$, and $S_1(\omega) > k\Delta_n$ and $|X_{(i+1)\Delta_n}(\omega) - X_{i\Delta_n}(\omega)| \leq 2\varepsilon$ for all $i \leq t/\Delta_n$ and $|X_{(i+1)k\Delta_n}(\omega) - X_{ik\Delta_n}(\omega)| \leq 2\varepsilon$ for all $i \leq t/k\Delta_n$. Next, on the set $\{(ik+j)\Delta_n < S_q \leq (ik+j+1)\Delta_n\}$ for $i \geq 1$ and $0 \leq j < k$, we set

$$L(n,q) = j,$$

$$K(n,q) = \frac{S_p}{\Delta_n} - (ik+j),$$

$$\alpha_-(n,q) = \frac{1}{\sqrt{\Delta_n}}(W_{S_q} - W_{(ik+j)\Delta_n}),$$

$$\alpha_+(n,q) = \frac{1}{\sqrt{\Delta_n}}(W_{(ik+j+1)\Delta_n} - W_{S_q}),$$

$$\beta_-(n,q) = \frac{1}{\sqrt{\Delta_n}}(W_{(ik+j)\Delta_n} - W_{ik\Delta_n}),$$

$$\beta_+(n,q) = \frac{1}{\sqrt{\Delta_n}}(W_{(i+1)k\Delta_n} - W_{(ik+j+1)\Delta_n}),$$

$$X(\varepsilon)_t = X_t - \sum_{q:\,S_q \leq t} \Delta X_{S_q},$$

$$R_-(n,q) = X(\varepsilon)_{(ik+j)\Delta_n} - X(\varepsilon)_{ik\Delta_n},$$

$$R'_+(n,q) = X(\varepsilon)_{(i+1)k\Delta_n} - X(\varepsilon)_{(ik+j+1)\Delta_n},$$

$$R_q^n = \Delta_{(ik+j+1)\Delta_n}^n X(\varepsilon), \qquad R_q'^n = R_q^n + R_-(n,q) + R_+(n,q),$$

$$R_q = \sqrt{\kappa_q} U_q \sigma_{S_q-} + \sqrt{1-\kappa_q} U_q' \sigma_{S_q},$$

$$R_-(q) = \sqrt{L_q}\,\overline{U}_q \sigma_{S_q-}, \qquad R_+(q) = \sqrt{k-1-L_q}\,\overline{U}_q' \sigma_{S_q}.$$

We can extend the proof of Lemma 6.2 of Jacod and Protter (1998) in this more complicated context, and with $\overset{\mathcal{L}-(s)}{\longrightarrow}$ denoting the stable convergence in law, to obtain that

$$(L(n,q), K(n,q), \alpha_-(n,q), \alpha_+(n,q), \beta_-(n,q), \beta_+(n,q))_{q \geq 1}$$

$$\overset{\mathcal{L}-(s)}{\longrightarrow} (L_q, \kappa_q, \sqrt{\kappa_q} U_q, \sqrt{1-\kappa_q} U_q', \sqrt{L_q}\,\overline{U}_q, \sqrt{k-L_q-1}\,\overline{U}_q')_{q \geq 1}.$$

We deduce from this, as in Lemma 5.9 of Jacod (2008), that

(55)
$$(R_q^n/\sqrt{\Delta_n}, R_-(n,q)/\sqrt{\Delta_n}, R_+(n,q)/\sqrt{\Delta_n})_{q \geq 1}$$
$$\overset{\mathcal{L}-(s)}{\longrightarrow} (R_q, R_-(q), R_+(q))_{q \geq 1}.$$

Now, since $f(x) = 0$ if $|x| \leq 2k\varepsilon$, on the set $\Omega_n(t,\varepsilon)$ and for all $s \leq t$ we have

$$V^n(f)_s - V(f)_{\Delta_n[s/\Delta_n]} = \sum_{q:\,S_q \leq \Delta_n[s/\Delta_n]} (f(\Delta X_{S_q} + R_q^n) - f(\Delta X_{S_q}))$$

$$= \sum_{q:\,S_q \leq \Delta_n[s/\Delta_n]} f'(\widetilde{R}_q^n) R_q^n,$$

where $\widetilde{R}_q^n$ is between $\Delta X_{S_q}$ and $\Delta X_{S_q} + R_q^n$. In a similar way, we get

$$\overline{V}^n(f)_s - V(f)_{k\Delta_n[s/k\Delta_n]} = \sum_{q:\,S_q \leq k\Delta_n[s/k\Delta_n]} f'(\widetilde{R}_q'^n) R_q'^n$$

with $\widetilde{R}'^n_q$ between $\Delta X_{S_q}$ and $\Delta X_{S_q} + R'^n_q$. Since $R^n_q$, $R_-(n,q)$ and $R_+(n,q)$ go to 0, we deduce that $\widetilde{R}^n_q$ and $\widetilde{R}'^n_q$ also go to 0 (for each $\omega$), whereas $\Omega_n(t,\varepsilon) \to \Omega$. Since $f'$ is continuous, we readily deduce the result of the theorem from (55).

Now we turn to the general case, where $f$ does not necessarily vanish around 0. With $\psi_\rho$ a $C^2$ function equal to 1 on $[-\rho, \rho]$ and to 0 outside $[-2\rho, 2\rho]$ and with $f_\rho = f\psi_\rho$, we have for all $\eta > 0$:

$$\lim_{\rho \to 0} \limsup_n \mathbb{P}\left(\sup_{s \leq t} |V^n(f_\rho)_s - V(f_\rho)_{\Delta_n[s/\Delta_n]}|/\sqrt{\Delta_n} > \eta\right) = 0$$

[this is (5.48) in Jacod (2008)] and we obviously have the same for $\overline{V}^n$. Since we have the stable convergence in law when we take the functions $f - f_\rho$, and since obviously $Z(f'_\rho)$ and $Z'(f'_\rho)$ converge locally uniformly in time to $Z(f')$ and $Z'(f')$, respectively, we readily deduce the stable convergence in law for the function $f$. □

We now prove part (a) of the theorem.

PROOF OF THEOREM 2(a) UNDER ASSUMPTION 2. The previous theorem applied with the function $f(x) = |x|^p$ (recall $p > 3$) yields the stable convergence in law of the processes

$$\left(\frac{1}{\sqrt{\Delta_n}}(\widehat{B}(p, \Delta_n)_t - B(p)_{\Delta_n[t/\Delta_n]}), \frac{1}{\sqrt{\Delta_n}}(\widehat{B}(p, k\Delta_n)_t - B(p)_{k\Delta_n[t/k\Delta_n]})\right)$$

toward $(Z(f')_t, Z(f')_t + Z'(f')_t)$, and $f'(x) = p|x|^{p-1}\text{sign}(x)$. Hence if $Z'(p, k)_t = Z'(f')_t$, the formula (15) follows from (53), and it remains to prove that $(B(p)_t - B(p)_{k\Delta_n[t/k\Delta_n]})/\sqrt{\Delta_n} \xrightarrow{\mathbb{P}} 0$ for all integers $k$. Since $|\delta(\omega, t, x)| \leq \gamma(x) \leq K$ we have

$$\mathbb{E}(B(p)_t - B(p)_{k\Delta_n[t/k\Delta_n]}) = \mathbb{E}\left(\sum_{k\Delta_n[t/k\Delta_n] < s \leq t} |\Delta X_s|^p\right)$$

$$= \mathbb{E}\left(\int_{k\Delta_n[t/k\Delta_n]}^t ds \int_E |\delta(s,x)|^p \lambda(dx)\right)$$

$$\leq \int_{k\Delta_n[t/k\Delta_n]}^t ds \int_E \gamma(x)^p \lambda(dx) \leq Kk\Delta_n,$$

where $K = \int \gamma(x)^p \lambda(dx)$ is finite (recall $p > 3$), and the result follows. □

For Theorem 2(b) it is also convenient to prove a more general result. Let $g = (g_j)_{1 \leq j \leq 2}$ be an $\mathbb{R}^2$-valued $C^2$ function on $\mathbb{R}^k$ with second partial derivatives having polynomial growth, and which is globally even [i.e.,

$g(-x) = g(x)$] and set

$$V^n(g)_t = \sum_{i=1}^{[t/k\Delta_n]} g(\Delta^n_{ki-k+1}X/\sqrt{\Delta_n}, \Delta^n_{ki-k+2}X/\sqrt{\Delta_n}, \ldots, \Delta^n_{ki}X/\sqrt{\Delta_n}). \quad (56)$$

We also denote by $\rho_y$ the normal law $N(0,y)$, and by $\rho_y^{k\otimes}$ its $k$-fold tensor product, and by $\rho_y^{k\otimes}(f)$ the integral of any function $f$ w.r.t. it. With the above assumptions we then have the following (a $d$-dimensional version is also available, of course):

THEOREM 9. *Under Assumption 1, assume in addition that $X$ is continuous. Then*

$$\Delta_n V^n(g)_t \xrightarrow{\mathbb{P}} V(g)_t := \frac{1}{k} \int_0^t \rho_{\sigma_u}^{k\otimes}(g). \quad (57)$$

*Furthermore the processes $\frac{1}{\sqrt{\Delta_n}}(\Delta_n V^n(g) - V(g))$ converge stably in law to the two-dimensional process $(Y(h,0), Y(h',h''))$ [recall (49)], where the processes $(h_t)$, $(h'_t)$ and $(h''_t)$ are such that the matrix $\begin{pmatrix} h_t & 0 \\ h'_t & h''_t \end{pmatrix}$ is a square root of the matrix $\Theta_t = (\Theta^{ij}_t)_{1 \leq i,j \leq 2}$ given by*

$$\Theta^{ij}_t = \frac{1}{k}(\rho_{\sigma_t}^{k\otimes}(g_i g_j) - \rho_{\sigma_t}^{k\otimes}(g_i)\rho_{\sigma_t}^{k\otimes}(g_j)). \quad (58)$$

PROOF. The proof is similar to the unipower case in Barndorff-Nielsen et al. (2006a). We indicate the changes that are necessary. Everywhere the sums over $i$ from 1 to $[nt]$ are replaced by sums from 1 to $[t/k\Delta_n]$. In (4.1) of that paper we replace $\beta^n_i = \sqrt{n}\sigma_{(i-1)/n}\Delta^n_i W$ by the collection $\beta^n_i(j) = \frac{1}{\sqrt{\Delta_n}}\sigma_{k(i-1)\Delta_n}\Delta^n_{k(i-1)+j}W$ for $j = 1, \ldots, k$. Then the same proof as for Proposition 4.1 of that paper shows that if

$$\widetilde{V}^n(g) = \sum_{i=1}^{[t/k\Delta_n]} g(\beta^n_i(1), \beta^n_i(2), \ldots, \beta^n_i(k)),$$

we have (57) and the stable convergence in law with $\widetilde{V}^n(g)$ instead of $V^n(g)$, and the same process $V(g)$. Next, similarly to Theorem 5.6 of that paper, and with $\zeta^n_i = g(\Delta^n_{ki-k+1}X/\sqrt{\Delta_n}, \ldots, \Delta^n_{ki}X/\sqrt{\Delta_n})$, the processes

$$\sqrt{\Delta_n} \sum_{i=1}^{[t/k\Delta_n]} (\zeta^n_i - \mathbb{E}(\zeta^n_i \mid \mathcal{F}_{k(i-1)\Delta_n}))$$

converge stably in law to $(Y(h,0), Y(h',h''))$, and this easily yields (57) for $V^n(g)$.

For the stable convergence in law, it remains to prove that the array

$$\zeta'^n_i = \sqrt{\Delta_n}\bigg(\mathbb{E}(\zeta^n_i \mid \mathcal{F}_{k(i-1)\Delta_n}) - \frac{1}{k\Delta_n}\int_{k(i-1)\Delta_n}^{ki\Delta_n} \rho^{k\otimes}_{\sigma_u}(g)\,du\bigg)$$

satisfies $\sum_{i=1}^{[t/\Delta_n]}|\zeta'^n_i| \xrightarrow{\mathbb{P}} 0$. This is proved as in Barndorff-Nielsen et al. (2006a); the fact that $g$ is a function of several variables makes no real difference, and since $g$ here is $C^2$ we are in the case of Hypothesis (K) of that paper [and not in the more complicated case of Hypothesis (K')]. □

Finally, we prove part (b) of the theorem.

PROOF OF THEOREM 2(b) UNDER ASSUMPTION 2. We simply apply the previous theorem to the even function $g$ with components

$$g_1(x_1,\ldots,x_k) = |x_1|^p + \cdots + |x_k|^p, \qquad g_2(x_1,\ldots,x_k) = |x_1 + \cdots + x_k|^p$$

(recall that $p \geq 2$, hence $g$ is $C^2$). The matrix $\Theta_t$ of (58) is then

$$\Sigma^{11}_t = (m_{2p} - m_p^2)|\sigma_t|^{2p}, \qquad \Sigma^{12}_t = (m_{k,p} - k^{p/2}m_p^2)|\sigma_t|^{2p},$$

$$\Sigma^{22}_t = k^{p-1}(m_{2p} - m_p^2)|\sigma_t|^{2p}.$$

Then a version of the triple $(h, h', h'')$ is given by $h_t = \alpha|\sigma_t|^p$ and $h'_t = \alpha'|\sigma_t|^p$ and $h''_t = \alpha''|\sigma_t|^p$, where

$$\alpha = \sqrt{m_{2r} - m_r^2}, \qquad \alpha' = \frac{1}{\alpha}(m_{k,r} - k^{r/2}m_r^2),$$

$$\alpha'' = \sqrt{k^{r-1}(m_{2r} - m_r^2) - \alpha\alpha'}.$$

Hence the result obtains, with $Y(p)_t = Y(h, 0)_t$ and $Y'(p, k)_t = Y(h', h'')_t$, upon noting that (16) follows from (50). □

8.3. *Proof of Theorem 3.*

PROOF. (a) Write $U_n = (\Delta_n)^{-1/2}(\widehat{B}(p,\Delta_n)_t - B(p)_t)$ and $V_n = (\Delta_n)^{-1/2} \times (\widehat{B}(p,k\Delta_n)_t - B(p)_t)$. Then

$$\widehat{S}(p,k,\Delta_n)_t - 1 = \frac{\widehat{B}(p,k\Delta_n)^n_t}{\widehat{B}(p,\Delta_n)_t} - 1 = (\Delta_n)^{1/2}\frac{V_n - U_n}{\widehat{B}(p,\Delta_n)_t}.$$

Theorem 3 yields that under Assumption 1 and if $p > 3$, then $V_n - U_n$ converges stably in law to $Z'(p,k)_t$, and the result follows from (15).

(b) Write $U'_n = (\Delta_n)^{-1/2}(\Delta_n^{1-p/2}\widehat{B}(p,\Delta_n) - m_p A(p))$ and $V'_n = (\Delta_n)^{-1/2} \times (\Delta_n^{1-p/2}\widehat{B}(p,k\Delta_n) - k^{p/2-1}m_p A(p))$. Then

$$\widehat{S}(p,k,\Delta_n)_t - k^{p/2-1} = \frac{\widehat{B}(p,k\Delta_n)^n_t}{\widehat{B}(p,\Delta_n)_t} - k^{p/2-1} = (\Delta_n)^{1/2}\frac{V'_n - k^{p/2-1}U'_n}{\widehat{B}(p,\Delta_n)_t}.$$

Since $V'_n - k^{p/2-1}U'_n$ converges stably in law to $Y'(p,k)_t - k^{p/2-1}Y(p)_t$ when $p \geq 2$ and $X$ is continuous, the result easily follows from (16). □

8.4. *Proof of Theorem 4.*

PROOF OF (25) UNDER ASSUMPTION 2. We fix $p \geq 2$, $t > 0$, $\varpi \in (1/2 - 1/p, 1/2)$ and $\alpha > 0$. Apply (47) to get, for any $B \geq 1$,

(59) $\qquad |U_n(B) - U'_n| \leq K_p(A_n(B) + A'_n(B) + A''_n(B))$,

where

$$U_n(B) = \Delta_n^{1-p/2} \sum_{i=1}^{[t/\Delta_n]} |\Delta_i^n X|^p 1_{\{|\Delta_i^n X| \leq \sqrt{B\Delta_n}\}}, \qquad U'_n = \Delta_n^{1-p/2} \sum_{i=1}^{[t/\Delta_n]} |\Delta_i^n X'|^p$$

and

$$A_n(B) = \Delta_n^{1-p/2} \sum_{i=1}^{[t/\Delta_n]} |\Delta_i^n X'|^p 1_{\{|\Delta_i^n X'| > \sqrt{B\Delta_n}/2\}},$$

$$A'_n(B) = B^{p/2-1} \sum_{i=1}^{[t/\Delta_n]} (|\Delta_i^n X''|^2 \wedge B\Delta_n),$$

$$A''_n(B) = \Delta_n^{1-p/2} \sum_{i=1}^{[t/\Delta_n]} |\Delta_i^n X'|^{p-1}(|\Delta_i^n X''| \wedge \sqrt{B\Delta_n}).$$

Then if we apply (44) and (45) with $\eta = \sqrt{B\Delta_n}$ and $\theta = \Delta_n^{1/4}$, and also Bienaymé–Chebyshev for $A_n(B)$ and Cauchy–Schwarz for $A''_n(B)$, and with the notation $\varepsilon_n^2 = \Delta_n^{1/2} + \Gamma(\Delta_n^{1/4})$ (so $\varepsilon_n \to 0$), we readily obtain (recall $B \geq 1$):

$$\mathbb{E}(A_n(B)) \leq \frac{K_p t}{B}, \qquad \mathbb{E}(A'_n(B)) \leq K_p t B^{p/2} \varepsilon_n^2, \qquad \mathbb{E}(A''_n(B)) \leq K_p t B^{1/2} \varepsilon_n.$$

It follows that (with other constants $K_p$):

$$\mathbb{P}(A_n(B) > B^{-1/2}) \leq K_p t B^{-1/2}, \qquad \mathbb{P}(A'_n(B) + A''_n(B) > \varepsilon_n^{1/2}) \leq K_p t \varepsilon_n^{1/2}.$$

On the other hand we can apply the last property in (11) to $X'$, which satisfies Assumption 2 and is continuous, to get $U'_n \xrightarrow{\mathbb{P}} m_p A(p)_t$. Combining this with the above estimates and (59), and since $\varepsilon_n \to 0$, we obtain for some constants $K$ and $K'$ depending on $p$ and $t$, and all $\eta > 0$ and $B > 1$:

(60) $\qquad \limsup_n \mathbb{P}(U_n(B) < m_p A(p)_t - KB^{-1/2} - \eta) \leq K'B^{-1/2}.$

Now, for any $B \geq 1$ we have $\alpha \Delta_n^\varpi > \sqrt{B\Delta_n}$ for all $n$ large enough because $\varpi < 1/2$. Therefore (60) remains valid if we replace the sets $\{|\Delta_i^n X| \leq$ $\sqrt{B\Delta_n}\}$ by $\{|\Delta_i^n X| \leq \alpha \Delta_n^{\varpi}\}$ in the definition of $U_n(B)$. Since $A(p)_t > 0$ a.s. and $B$ is arbitrarily large and $\eta$ arbitrarily small, we deduce that $C_n = \sum_{i=1}^{[t/\Delta_n]} |\Delta_i^n X|^p 1_{\{|\Delta_i^n X| \leq \alpha \Delta_n^{\varpi}\}}$ satisfies

$$(61) \quad \mathbb{P}\bigg(\frac{1}{\Delta_n^{1-p/2} C_n} > \frac{2}{m_p A(p)_t}\bigg) = \mathbb{P}\bigg(\Delta_n^{1-p/2} C_n < \frac{m_p A(p)_t}{2}\bigg) \to 0.$$

At this stage, the proof of (25) is straightforward. Indeed, since $|\Delta_i^n X|^{2p} \leq \alpha^p \Delta_n^{p\varpi} |\Delta_i^n X|^p$ when $|\Delta_i^n X| \leq \alpha \Delta^{\varpi}$, one deduces from (21) that

$$\frac{\Delta_n \widehat{A}(2p, \Delta_n)_t}{\widehat{A}(p, \Delta_n)_t^2} \leq \frac{K \Delta_n^{p\varpi}}{C_n} = \frac{K \Delta_n^{p\varpi+1-p/2}}{\Delta_n^{1-p/2} C_n}.$$

Since $p\varpi + 1 - p/2 > 0$, the result readily follows from (61). $\square$

PROOF OF (26) UNDER ASSUMPTION 2. Barndorff-Nielsen, Shephard and Winkel (2006b) proved a similar result under more restrictive conditions. Apply (45) with $\eta = \sqrt{\Delta_n}$ and $\theta = \Delta_n^{1/4}$, and divide par $\Delta_n$, to get

$$(62) \quad \mathbb{E}_{i-1}^n(|\Delta_i^n X''/\sqrt{\Delta_n}|^2 \wedge 1) \leq \alpha_n$$

for a (deterministic) sequence $\alpha_n$ going to 0. Next, we apply (48) with $0 < r < 2$ and $x = \Delta_i^n X'/\sqrt{\Delta_n}$ and $y = \Delta_i^n X''/\sqrt{\Delta_n}$ and we use (44) for $q = 2$ and (62) to get

$$\mathbb{E}_{i-1}^n(||\Delta_i^n X/\sqrt{\Delta_n}|^r - |\Delta_i^n X'/\sqrt{\Delta_n}|^r|) \leq K(\varepsilon^r + A\varepsilon + A^{r-2} + A^r \varepsilon^{-2} \alpha_n).$$

Hence by taking $\varepsilon = \varepsilon_n = \alpha_n^{1/4}$ and $A = A_n = \alpha_n^{-1/8}$, we get as soon as $\alpha_n \leq 1$:

$$(63) \quad \Delta_n^{-r/2} \sup_i \mathbb{E}_{i-1}^n(||\Delta_i^n X|^r - |\Delta_i^n X'|^r|) \leq \alpha_n' = K(\alpha_n^{1/8} + \alpha^{(2-r)/8}) \to 0.$$

Now, from Barndorff-Nielsen et al. (2006a), we know that the result holds when $X$ is continuous: the processes $V'(r, q, \Delta_n)$ defined by (22) with $X$ substituted with $X'$ converge in probability, locally uniformly in time, to the process $A(qr)$ (this holds even when $r \geq 2$). Therefore it is enough to prove that

$$(64) \quad \mathbb{E}\bigg(\sup_{s \leq t} |V(q, r, \Delta_n)_s - V'(q, r, \Delta_n)_s|\bigg) \to 0.$$

The left-hand side of (64) at time $t$ is smaller than $\sum_{i=1}^{[t/\Delta_n]} E(\zeta_i^n)$, where

$$\zeta_i^n = \frac{\Delta_n^{1-qr/2}}{m_r^q}\bigg|\prod_{j=1}^q |\Delta_{i+j-1}^n X|^r - \prod_{j=1}^{q-1} |\Delta_{i+j-1}^n X'|^r\bigg| = \frac{1}{m_r^q} \sum_{l=1}^q \zeta_i^n(l),$$

$$\zeta_i^n(l) = \Delta_n^{1-rq/2} \prod_{j=1}^{l-1} |\Delta_{i+j-1}^n X|^r ||\Delta_{i+l-1}^n X|^r - |\Delta_{i+l-1}^n X'|^r| \prod_{j=l+1}^q |\Delta_{i+j-1}^n X'|^r,$$

where an empty product is set to be 1. Taking successive conditional expectations, and using (44) and (63), we readily obtain that $E(\zeta_i^n(l)) \leq K\alpha_n' \Delta_n$. Then obviously $\sum_{i=1}^{[t/\Delta_n]} E(\zeta_i^n(l)) \leq Kt\alpha_n' \to 0$ for each $l = 1, \ldots, q$, hence (64). □

LEMMA 1. *Under Assumption 2 and if*

$$(65) \qquad \widehat{D}_t'^n := \frac{1}{k_n \Delta_n} \sum_{i=1}^{[t/\Delta_n]} |\Delta_i^n X|^p \sum_{j \in I_{n,t}(i)} (\delta_j^n)^2 \xrightarrow{\mathbb{P}} D(p)_t$$

*and also, when $X$ is continuous,*

$$(66) \qquad \widehat{D}_t''^n := \frac{1}{k_n \Delta_n^{p/2}} \sum_{i=1}^{[t/\Delta_n]} |\delta_i^n|^p \sum_{j \in I_{n,t}(i)} (\delta_j^n)^2 \xrightarrow{\mathbb{P}} m_p A(p+2)_t,$$

*then (27) holds, as well as (28) when $X$ is continuous.*

PROOF. In (28) both members are increasing in $t$ and the limit is continuous in $t$, so it is enough to prove the convergence separately for each $t$. To unify the proof of the two results, we write $u_n = 1$ and $D = D(p)$ if $X$ jumps, and $u_n = \Delta_n^{1-p/2}$ and $D = m_p A(p+2)$ if $X$ is continuous, in which case we also set

$$\widehat{D}_t'^n := \frac{u_n}{k_n \Delta_n} \sum_{i=1}^{[t/\Delta_n]} |\Delta_i^n X|^p \sum_{j \in I_{n,t}(i)} (\delta_j^n)^2.$$

With this notation, (27) and (28) amount to $u_n \widehat{D}(p, \Delta_n)_t \xrightarrow{\mathbb{P}} D_t$, and in a first step we show that this is equivalent to $\widehat{D}_t'^n \xrightarrow{\mathbb{P}} D_t$.

Apply (47) with $p = 2$ and $x = \delta_i^n$ and $y = \Delta_i^n X - \delta_i^n$ and $\theta = \alpha \Delta_n^\varpi$, and use Cauchy–Schwarz and Bienaymé–Chebyshev and (44), and also (46) with $\eta = \eta_n = \alpha \Delta_n^\varpi$ and $\varepsilon = \eta_n^{1/2}$, to get after some simple calculations that, with $K = K_\alpha$:

$$(67) \quad \begin{aligned} &\frac{1}{\Delta_n} \mathbb{E}(|(\theta_i^n)^2 - (\delta_i^n)^2| \mid \mathcal{F}_{(i-1)\Delta_n}) \\ &\qquad \leq \widehat{\Gamma}_n := K(\Delta_n + \Delta_n^\varpi + \Gamma(\eta_n^{1/2}) + \Delta_n^{\varpi/2} + \sqrt{\Gamma(\eta_n^{1/2})}) \to 0, \end{aligned}$$

where the final convergence follows from $\eta_n \to 0$, hence $\Gamma(\eta_n^{1/2}) \to 0$ as well. Observe that $u_n \widehat{D}(p, \Delta_n)$ is the same as $\widehat{D}'^n$, with $\delta_j^n$ substituted with $\theta_i^n$. Hence the difference $u_n \widehat{D}(p, \Delta_n)_t - \widehat{D}_t'^n$ is the sum of less than $2k_n[t/\Delta_n]$ terms (strictly less, because of the border effects at 0 and $[t/\Delta_n]$), each one being smaller than $\frac{u_n}{k_n\Delta_n}|\Delta_i^n X|^p|(\theta_j^n)^2 - (\delta_j^n)^2|$, for some $i \neq j$. Then, by taking two successive conditional expectations and using (44) and (67), we see that the expectation of such a term is smaller than $K_p\Delta_n\widehat{\Gamma}_n/k_n$ in both cases ($X$ continuous or not). Therefore $E(|u_n\widehat{D}(p,\Delta_n)_t - \widehat{D}'^n_t|) \leq K_p t\widehat{\Gamma}_n$. Since $\widehat{\Gamma}_n \to 0$, the claim of this step is complete.

So far we have proved that (65) implies (27). For proving that (66) implies (28) when $X$ is continuous, it remains to show that in this case $E(|\widehat{D}'^n_t - \widehat{D}''^n_t|) \to 0$. Exactly as above, $\widehat{D}'^n_t - \widehat{D}''^n_t$ is the sum of less than $2k_n[t/\Delta_n]$ terms, each one being smaller than $\zeta'_{i,j} = \frac{u_n}{k_n\Delta_n}||\Delta_i^n X|^p - |\delta_j^n|^p|(\delta_j^n)^2$, for some $i \neq j$. Since $||x+y|^p - |x|^p| \leq K_p(|x|^{p-1}|y| + |y|^p)$, and since $X = X'$ because $X$ is continuous, we deduce from (44) and by taking two successive conditional expectations that $E(\zeta_{i,j}^n) \leq K_p u_n \Delta_n^{1/2+p/2}/k_n$, hence $E(|\widehat{D}'^n_t - \widehat{D}''^n_t|) \leq K_p t\sqrt{\Delta_n}$, and we are done. □

PROOF OF (27) UNDER ASSUMPTION 2. *Step* 1. For any $\rho \in (0,1)$ we set $\psi_\rho(x) = 1 \wedge (2-|x|/\rho)^+$ (a continuous function with values in $[0,1]$, equal to $1$ if $|x| \leq \rho$ and to $0$ if $|x| \geq 2\rho$), and we introduce two increasing processes:

$$Y(\rho)_t^n = \frac{1}{k_n\Delta_n} \sum_{i=1}^{[t/\Delta_n]} \psi_\rho(\Delta_i^n X)|\Delta_i^n X|^p \sum_{j \in I_{n,t}(i)} (\delta_j^n)^2, \qquad Z(\rho)_t^n = \widehat{D}'^n_t - Y(\rho)_t^n.$$

By the previous lemma we need to prove (65), and for this it is obviously enough to show the following three properties, for some suitable processes $Z(\rho)$:

(68) $$\lim_{\rho \to 0} \limsup_n \mathbb{E}(Y(\rho)_t^n) = 0,$$

(69) $$\rho \in (0,1), \ n \to \infty \quad \Rightarrow \quad Z(\rho)_t^n \xrightarrow{\mathbb{P}} Z(\rho)_t,$$

(70) $$\rho \to 0 \quad \Rightarrow \quad Z(\rho)_t \xrightarrow{\mathbb{P}} D_t.$$

*Step* 2. By singling out the cases $2|x| > |y|$ and $2|x| \leq |y|$, we check that (recall $\rho < 1$):

$$\psi_\rho(x+y)|x+y|^p \leq K_p(|x|^p + (y^2 \wedge \rho^2)).$$

Using this with $x = \delta_i^n$ and $y = \Delta_i^n X - \delta_i^n$, plus (44) and (46) with $\eta = \rho$ and $\theta = \sqrt{\rho}$, we obtain

$$\frac{1}{\Delta_n}\mathbb{E}_{i-1}^n(\psi_\rho(\Delta_i^n X)|\Delta_i^n X|^p) \leq \Gamma'_n(\rho) := K_p(\Delta_n^{p/2-1} + \rho + \Gamma(\sqrt{\rho}) + \Delta_n\rho^{-1}).$$

Now, $Y(\rho)_t^n$ is the sum of less than $2k_n[t/\Delta_n]$ terms, all smaller than $\frac{1}{k_n\Delta_n} \times |\Delta_i^n X|^p|\delta_j^n|^2$ for some $i \neq j$. By taking two successive conditional expectations, as in the previous lemma, and by using (44) and the above, we see

that the expectation of such a term is smaller than $K_p \Delta_n \Gamma'_n(\rho)/k_n$. Thus $E(Y(\rho)^n_t) \leq Kt\Gamma'_n(\rho)$, and since obviously $\lim_{\rho \to 0} \limsup_n \Gamma'_n(\rho) = 0$ we obtain (68).

*Step* 3. Now we define the process $Z(\rho)$. Let us call $T_q(\rho)$ for $q = 1, 2, \ldots$ the successive jump times of the Poisson process $\mu([0, t] \times \{x : \gamma(x) > \rho/2\})$, and set

$$Z(\rho)_t = \sum_{q \,:\, T_q(\rho) \leq t} |\Delta X_{T_q(\rho)}|^p (1 - \psi_\rho(\Delta X_{T_p(\rho)}))(\sigma^2_{T_q(\rho)-} + \sigma^2_{T_q(\rho)}).$$

For all $\omega \in \Omega$, $q \geq 1$, $\rho' \in (0, \rho)$ there is $q'$ such that $T_q(\rho)(\omega) = T_{q'}(\rho')(\omega)$, whereas $1 - \psi_\rho$ increases to the indicator of $\mathbb{R} \setminus \{0\}$. Thus $Z(\rho)_t(\omega) \uparrow D(p)_t(\omega)$, and we have (70).

*Step* 4. It remains to prove (69). Fix $\rho \in (0, 1)$ and write $T_q = T_q(\rho)$. Recall that for $u$ different from all $T_q$'s, we have $|\Delta X_u| \leq \rho/2$. Hence, for each $\omega$ and each $t > 0$, we have the following properties for all $n$ large enough: there is no $T_q$ in $(0, k_n\Delta_n]$, nor in $(t - (k_n + 1)\Delta_n, t]$; there is at most one $T_q$ in an interval $((i-1)\Delta_n, i\Delta_n]$ with $i\Delta_n \leq t$, and if this is not the case we have $\psi_\rho(\Delta^n_i X) = 1$. Hence for $n$ large enough we have

$$Z(\rho)^n_t = \sum_{q \,:\, k_n\Delta_n < T_q \leq t - (k_n+1)\Delta_n} \zeta^n_q,$$

where

$$\zeta^n_q = \frac{1}{k_n\Delta_n} |\Delta^n_{i(n,q)} X|^p (1 - \psi_\rho(\Delta^n_{i(n,q)} X)) \sum_{j \in I'(n,q)} (\delta^n_j)^2$$

and $i(n, q) = \inf(i : i\Delta_n \geq T_q)$ and $I'(n, q) = \{j : j \neq i(n, q), |j - i(n, q)| \leq k_n\}$.

To get (69) it is enough that $\zeta^n_q \xrightarrow{\mathbb{P}} |\Delta X_{T_q}|^p (1 - \psi_\rho(\Delta X_{T_q}))(\sigma^2_{T_q-} + \sigma^2_{T_q})$ for any $q$. Since $|\Delta^n_{i(n,q)} X|^p (1 - \psi_\rho(\Delta^n_{i(n,q)} X)) \to |\Delta X_{T_q}|^p (1 - \psi_\rho(\Delta X_{T_p}))$ pointwise, so it remains to prove that

$$(71) \quad \frac{1}{k_n\Delta_n} \sum_{j \in I'_-(n,q)} (\delta^n_j)^2 \xrightarrow{\mathbb{P}} \sigma^2_{T_q-}, \qquad \frac{1}{k_n\Delta_n} \sum_{j \in I'_+(n,q)} (\delta^n_j)^2 \xrightarrow{\mathbb{P}} \sigma^2_{T_q},$$

where $I'_-(n, q)$ and $I'_+(n, q)$ are the subsets of $I'(n, q)$ consisting in those $j$ smaller, respectively bigger, than $i(n, q)$. We write

$$U^n_q = \frac{1}{k_n\Delta_n} \sum_{j \in I'_-(n,q)} (\Delta^n_j W)^2,$$

$$s^n_q = \inf_{u \in [T_q - k_n\Delta_n, T_q)} \sigma^2_u,$$

$$S^n_q = \sup_{u \in [T_q - k_n\Delta_n, T_q)} \sigma^2_u.$$

On the one hand, both $s_q^n$ and $S_q^n$ converge as $n \to \infty$ to $\sigma_{T_q-}^2$, because $k_n \Delta_n \to 0$. On the other hand, the left-hand side of the first expression in (71) is in between the two quantities $s_q^n U_q^n$ and $S_q^n U_q^n$. Moreover the variables $\Delta_i^n W$ are i.i.d. $N(0, \Delta_n)$, so $U_q^n$ is distributed as $U_n' = \frac{1}{k_n} \sum_{i=1}^{k_n} V_i$ where the $V_i$'s are i.i.d. $N(0,1)$, and by the usual law of large numbers we have $U_n' \to 1$ a.s., hence $U_q^n \xrightarrow{\mathbb{P}} 1$; these facts entail the first part of (71), and the second part is proved in the same way. $\square$

PROOF OF (28) UNDER ASSUMPTION 2. *Step* 1. We have to prove that when $X$ is continuous, then (66) holds. We have

$$\Delta_n^{-p/2} \sum_{i=1}^{[t/\Delta_n]} |\delta_i^n|^{p+2} \xrightarrow{\mathbb{P}} m_{p+2} A(p+2)_t$$

[this property is the preliminary step in Jacod (2008) to prove the last part of (11); alternatively, it can be deduced from (11) exactly as $\widehat{D}_t'^n - \widehat{D}_t''^n \xrightarrow{\mathbb{P}} 0$ in the end of the proof of Lemma 1]. Therefore it is enough to prove that

$$Y_t^n := \frac{1}{k_n \Delta_n^{p/2}} \sum_{i=1}^{[t/\Delta_n]} \sum_{j \in I_{n,t}(i)} \zeta_{i,j}^n \xrightarrow{\mathbb{P}} 0,$$

where $\zeta_{i,j}^n = m_{p+2} |\delta_i^n|^p |\delta_j^n|^2 - m_p |\delta_i^n|^{p+2}$.

Let $T(n,i) = (i - k_n - 1)^+ \Delta_n$. When $j > i$ we have (recall $\sigma_t$ is bounded):

$$|\mathbb{E}_{i-1}^n(\zeta_{i,j}^n)| = |\mathbb{E}_{i-1}^n(m_{p+2} |\sigma_{(i-1)\Delta_n}|^p (\sigma_{(j-1)\Delta_n}^2 - \sigma_{(i-1)\Delta_n}^2) |\Delta_i^n W|^p |\Delta_j^n W|^2$$
$$+ |\sigma_{(i-1)\Delta_n}|^{p+2} (m_{p+2} |\Delta_i^n W|^p |\Delta_j^n W|^2 - m_p |\Delta_i^n W|^{p+2}))|$$
$$\leq K_p \Delta_n \mathbb{E}_{i-1}^n (|\sigma_{(j-1)\Delta_n}^2 - \sigma_{(i-1)\Delta_n}^2 ||\Delta_i^n W|^p|).$$

When $i - k_n \leq j < i$ we have

$$|\mathbb{E}(\zeta_{i,j}^n \mid \mathcal{F}_{(j-1)\Delta_n})|$$
$$= |\mathbb{E}(m_{p+2} |\sigma_{(j-1)\Delta_n}|^2 (|\sigma_{(i-1)\Delta_n}|^p - |\sigma_{(j-1)\Delta_n}|^p) |\Delta_i^n W|^p |\Delta_j^n W|^2$$
$$+ |\sigma_{(j-1)\Delta_n}|^{p+2} (m_{p+2} |\Delta_i^n W|^p |\Delta_j^n W|^2 - m_p |\Delta_i^n W|^{p+2}) \mid \mathcal{F}_{(j-1)\Delta_n})|$$
$$\leq K_p \Delta_n^{p/2} \mathbb{E}(||\sigma_{(i-1)\Delta_n}|^p - |\sigma_{(j-1)\Delta_n}|^p| |\Delta_j^n W|^2 |\mathcal{F}_{(j-1)\Delta_n}).$$

Therefore, since $||\sigma_{(i-1)\Delta_n}|^q - |\sigma_{(j-1)\Delta_n}|^q| \leq K_q |\sigma_{(i-1)\Delta_n} - \sigma_{(j-1)\Delta_n}|$ for any $q \geq 1$, we deduce from (44) and Cauchy–Schwarz and the two estimates above that

(72) $\qquad j \in I_{n,t}(i) \quad \Rightarrow \quad |\mathbb{E}(\zeta_{i,j}^n \mid \mathcal{F}_{T(n,i)})| \leq K_p \Delta_n^{3/2+p/2}.$

Moreover, as a trivial consequence of (44) again, we get

(73) $$j \in I_{n,t}(i) \quad \Rightarrow \quad \mathbb{E}(|\zeta_{i,j}^n|^2 \mid \mathcal{F}_{T(n,i)}) \leq K_p \Delta_n^{2+p}.$$

Now we set

$$\eta_i^n = \frac{1}{k_n \Delta_n^{p/2}} \sum_{j \in I_{n,t}} \zeta_{i,j}^n, \qquad Z_n = \sum_{i=1}^{[t/\Delta_n]} \mathbb{E}(\eta_i^n \mid \mathcal{F}_{T(n,i)}), \qquad Z_n' = Y_t^n - Z_n.$$

On the one hand (72) yields $|Z_n| \leq K_p t \sqrt{\Delta_n} \to 0$. On the other hand $Z_n' = \sum_{i=1}^{[t/\Delta_n]} (\eta_i^n - E(\eta_i^n \mid \mathcal{F}_{T(n,i)}))$, whereas $\eta_i^n$ is $\mathcal{F}_{T(n,i+2k_n+1)}$-measurable. Hence (73) gives

$$\mathbb{E}(Z_n'^2) = \sum_{i,i' : 1 \leq i, i' \leq [t/\Delta_n]} \mathbb{E}((\eta_i^n - \mathbb{E}(\eta_i^n \mid \mathcal{F}_{T(n,i)}))(\eta_{i'}^n - \mathbb{E}(\eta_{i'}^n \mid \mathcal{F}_{T(n,i')})))$$

$$\leq \sum_{i,i' : 1 \leq i \leq [t/\Delta_n], |i-i'| \leq 2k_n+1} \mathbb{E}((\eta_i^n - \mathbb{E}(\eta_i^n \mid \mathcal{F}_{T(n,i)}))$$

$$\times (\eta_{i'}^n - \mathbb{E}(\eta_{i'}^n \mid \mathcal{F}_{T(n,i')})))$$

$$\leq \sum_{i,i' : 1 \leq i \leq [t/\Delta_n], |i-i'| \leq 2k_n+1} \mathbb{E}(|\eta_i^n \eta_{i'}^n|) \leq (4k_n + 3)[t/\Delta_n] K_p \Delta_n^2$$

$$\leq K_p t (k_n \Delta_n),$$

which goes to 0 by virtue of (23). It then follows that $Y_t^n \xrightarrow{\mathbb{P}} 0$, and we are done. $\square$

### 8.5. Proof of Theorem 6.

PROOF. (a) When $X$ is continuous the variables $U_n = (\widehat{V}_{n,t}^c)^{-1/2}(\widehat{S}(p,k,\Delta_n)_t - k^{p/2-1})$ converge stably in law to $N(0,1)$ by Theorem 5(b), for both choices of $\widehat{V}_{n,t}^c$. In particular if $c_{n,t}^c$ is given by (35) we have that $\mathbb{P}(C_{n,t}^c) = \mathbb{P}(U_n < -z_\alpha) \to \alpha$, and we even have more, because of the stable convergence in law; namely, for any measurable subset $B$ of $\Omega$, then

(74) $$\mathbb{P}(C_{n,t}^c \cap B) = \mathbb{P}(\{U_n < -z_\alpha\} \cap B) \to \alpha \mathbb{P}(B).$$

Now, this is not quite enough to prove (a), since it may happen that $X$ is not continuous but the observed path is continuous on $[0, t]$. To deal with this case, for any integer $N$ we introduce the process

$$X_t^{(N)} = X_0 + \int_0^t b_s \, ds + \int_0^t \sigma_s \, dW_s - \int_0^t \delta_s' 1_{\{|\delta_s'| \leq N\}} \, ds.$$

Put an additional exponent $^{(N)}$ for the variables defined on the basis of $X^{(N)}$, writing, for example, $\widehat{V}_{n,t}^{c(N)}$ or $\widehat{S}(p,k,\Delta_n)_t^{(N)}$ or $C_{n,t}^{c(N)}$. Then (74) applied with the continuous process $X^{(N)}$ shows that $\mathbb{P}(C_{n,t}^{c(N)} \cap B) \to \alpha \mathbb{P}(B)$ for any $B \in \mathcal{F}$. However, on the set $\Omega_t^{(N)}$ where $X_s = X_s^{(N)}$ for all $s \in [0, t]$ we have $\widehat{V}_{n,t}^{c(N)} = \widehat{V}_{n,t}^c$ and $\widehat{S}(p, k, \Delta_n)_t^{(N)} = \widehat{S}(p, k, \Delta_n)_t$, hence $C_{n,t}^{c(N)} \cap \Omega_t^{(N)} = C_{n,t}^c \cap \Omega_t^{(N)}$. Therefore $\mathbb{P}(C_{n,t}^c \cap \Omega_t^{(N)}) \to \alpha \mathbb{P}(\Omega_t^{(N)})$. Since $\Omega_t^{(N)}$ increases to $\Omega_t^c$ by virtue of (e) of Assumption 1, we deduce (i).

(b) Finally assume $\mathbb{P}(\Omega_t^j) > 0$. Then Theorem 1 yields that $\widehat{S}(p, k, \Delta_n)_t \xrightarrow{\mathbb{P}} 1$ on the set $\Omega_t^j$. On the other hand if we use the version (30) for $\widehat{V}_{n,t}^c$ we deduce from (25) that $c_{n,t}^c \xrightarrow{\mathbb{P}} k^{p/2-1} > 1$, whereas if we use the version (31) we have $\widehat{V}_{n,t}^c/\Delta_n \xrightarrow{\mathbb{P}} M(p, k)A(2p)_t/A(p)_t^2$, hence again $c_{n,t}^c \xrightarrow{\mathbb{P}} k^{p/2-1} > 1$; so the result is obvious. $\square$

### 8.6. Proof of Theorem 7.

PROOF. Relative to the proof of the previous theorem, we need a few changes. First we replace (a) of Theorem 6 by two statements (a) and (b) here: the case (b) corresponds to the situation where the limit in Theorem 5(a) is normal, and so this is similar to Theorem 6(a). The case (a) here corresponds to a nonnormal limit with variance 1, and for this limit we cannot evaluate exactly the quantiles and we rely upon the Chebyshev inequality; this is why we only get a bound on the level but not the exact value. Apart from these changes, the proof for the level is the same.

For (c) we suppose that $\mathbb{P}(\Omega_t^c) > 0$. Letting $T = \inf(s : \Delta X_s \neq 0)$, we have $T > t$ on $\Omega_t$, whereas $\mu(\{(s, z) : s \leq T, \delta(s, z) \neq 0\}) \leq 1$, hence $\mathbb{E}(\nu(\{(s, z) : s \leq T, \delta(s, z) \neq 0\})) \leq 1$, hence a fortiori the predictable process $Y_s = \int_0^s \int_E |\kappa \circ \delta(r, z)|\nu(dr, dz)$ is finite-valued on $[0, T]$, so there is an increasing sequence $(T_q)$ of stopping times with $T_q \leq T$ and $Y_{T_q} \leq q$ and $\mathbb{P}(\{T_q < t\} \cap \Omega_{n,t}^c) \to \mathbb{P}(\Omega_{n,t}^c)$, and it is clearly enough to show that $\mathbb{P}(C_{n,t}^j | \{T_q < t\} \cap \Omega_{n,t}^c) \to 1$ for all $q$ having $\mathbb{P}(\{T_q < t\} \cap \Omega_{n,t}^c) > 0$.

Now with $q$ fixed we consider the process

$$\overline{X}_s = X_0 + \int_0^s b_r\, dr + \int_0^s \sigma_r\, dW_r - \int_0^{s \wedge T_q} \int_E \kappa \circ \delta(r, z)\nu(dr, dz).$$

This process satisfies Assumption 1 and is continuous, and it coincides with $X$ on the interval $[0, T_q]$. Then similarly to the previous proof, we see that on the set $\{T_q < t\} \cap \Omega_{n,t}^c$ we have $\widehat{S}(p, k, \Delta_n)_t \xrightarrow{\mathbb{P}} k^{p/2-1}$ and $\widehat{V}_{n_t}^j \xrightarrow{\mathbb{P}} 0$ [use (23) and (28) for the latter], hence the result. $\square$

**Acknowledgments.** We are very grateful to a referee, an Associate Editor and Cecilia Mancini, for constructive comments that helped us greatly improve this paper.

## REFERENCES

AÏT-SAHALIA, Y. (2002). Telling from discrete data whether the underlying continuous-time model is a diffusion. *J. Finance* **57** 2075–2112.

AÏT-SAHALIA, Y. (2004). Disentangling diffusion from jumps. *J. Financial Economics* **74** 487–528.

AÏT-SAHALIA, Y. and JACOD, J. (2007). Volatility estimators for discretely sampled Lévy processes. *Ann. Statist.* **35** 355–392. MR2332279

AÏT-SAHALIA, Y. and KIMMEL, R. (2007). Maximum likelihood estimation of stochastic volatility models. *J. Financial Economics* **83** 413–452.

ANDERSEN, T. G., BOLLERSLEV, T. and DIEBOLD, F. X. (2003). Some like it smooth, and some like it rough. Technical Report, Northwestern Univ.

BARNDORFF-NIELSEN, O., GRAVERSEN, S., JACOD, J., PODOLSKIJ, M. and SHEPHARD, N. (2006a). A central limit theorem for realised bipower variations of continuous semimartingales. In *From Stochastic Calculus to Mathematical Finance, The Shiryaev Festschrift* (Y. Kabanov, R. Liptser and J. Stoyanov, eds.) 33–69. Springer, Berlin. MR2233534

BARNDORFF-NIELSEN, O. E. and SHEPHARD, N. (2006). Econometrics of testing for jumps in financial economics using bipower variation. *J. Financial Econometrics* **4** 1–30.

BARNDORFF-NIELSEN, O. E., SHEPHARD, N. and WINKEL, M. (2006b). Limit theorems for multipower variation in the presence of jumps. *Stochastic Process. Appl.* **116** 796–806. MR2218336

CARR, P. and WU, L. (2003). What type of process underlies options? A simple robust test. *J. Finance* **58** 2581–2610.

HUANG, X. and TAUCHEN, G. (2006). The relative contribution of jumps to total price variance. *J. Financial Econometrics* **4** 456–499.

JACOD, J. (2008). Asymptotic properties of realized power variations and related functionals of semimartingales. *Stochastic Process. Appl.* **118** 517–559.

JACOD, J. and PROTTER, P. (1998). Asymptotic error distributions for the Euler method for stochastic differential equations. *Ann. Probab.* **26** 267–307. MR1617049

JACOD, J. and SHIRYAEV, A. N. (2003). *Limit Theorems for Stochastic Processes*, 2nd ed. Springer, New York. MR1943877

JIANG, G. J. and OOMEN, R. C. (2005). A new test for jumps in asset prices. Technical report, Univ. Warwick, Warwick Business School.

LEE, S. and MYKLAND, P. A. (2008). Jumps in financial markets: A new nonparametric test and jump clustering. Review of Financial Studies. To appear.

LEPINGLE, D. (1976). La variation d'ordre p des semi-martingales. *Z. Wahrsch. Verw. Gebiete* **36** 295–316. MR0420837

MANCINI, C. (2001). Disentangling the jumps of the diffusion in a geometric jumping Brownian motion. *Giornale dell'Istituto Italiano degli Attuari* **LXIV** 19–47.

MANCINI, C. (2004). Estimating the integrated volatility in stochastic volatility models with Lévy type jumps. Technical report, Univ. Firenze.

RÉNYI, A. (1963). On stable sequences of events. *Sankyā Ser. A* **25** 293–302. MR0170385

WOERNER, J. H. (2006a). Analyzing the fine structure of continuous-time stochastic processes. Technical Report, Univ. Göttingen.

WOERNER, J. H. (2006b). Power and multipower variation: Inference for high-frequency data. In *Stochastic Finance* (A. Shiryaev, M. do Rosário Grosshino, P. Oliviera and M. Esquivel, eds.) 264–276. Springer, Berlin.

Department of Economics  
Princeton University and NBER  
Princeton, New Jersey 08544-1021  
USA  
E-mail: yacine@princeton.edu

Institut de Mathématiques de Jussieu  
CNRS UMR 7586  
Université Pierre et Marie Curie  
75252 Paris Cédex 05  
France  
E-mail: jean.jacod@upmc.fr